\documentclass[11pt]{article}
\usepackage{amssymb,latexsym,times}
\usepackage{mathptmx}
\usepackage{amscd}
\usepackage{MnSymbol}
\usepackage[all]{xy}

\catcode`\@=11
\renewcommand\subsection{\@startsection{subsection}{2}{\z@}%
                                     {-3.25ex\@plus -1ex \@minus -.2ex}%
                                     {-0.01 mm}
                                     {\normalfont\large\bfseries}}

\catcode`\@=11
\renewcommand\subsubsection{\@startsection{subsubsection}{2}{\z@}%
                                     {-3.25ex\@plus -1ex \@minus -.2ex}%
                                     {-0.01 mm}
                                     {\normalfont\bfseries}}

\newtheorem{Thm}{Theorem}[section]
\newtheorem{Lem}[Thm]{Lemma}
\newtheorem{Cor}[Thm]{Corollary}
\newtheorem{Prop}[Thm]{Proposition}

\newtheorem{example}[Thm]{Example}
\newtheorem{remark}[Thm]{Remark}
\newtheorem{Def}[Thm]{Definition}

\newcommand{\A}{\mathcal{A}}

\newcommand{\E}{\mathcal{E}}
\newcommand{\Z}{\mathbb{Z}}
\newcommand{\Q}{\mathbb{Q}}
\newcommand{\RR}{\mathcal{R}}

\newcommand{\N}{\mathbb{N}}
\newcommand{\C}{\mathbb{C}}

\newcommand{\md}{\operatorname{mod}}

\newcommand{\Hom}{\operatorname{Hom}} 
\newcommand{\Ext}{\operatorname{Ext}}
\newcommand{\End}{\operatorname{End}}

\newcommand{\bsm}{\begin{smallmatrix}}
\newcommand{\esm}{\end{smallmatrix}}

\newcommand{\bd}{\mathbf{d}}

\newcommand{\Cok}{\operatorname{Coker}}
\def\Sub{\mathrm{Sub}}
\def\Fac{\mathrm{Fac}}

\def\proof{\medskip\noindent {\it Proof --- \ }}
\def\cqfd{\hfill $\Box$ \bigskip}
\def\CC{{\mathcal C}}
\def\resp{{\em resp.\ }}

\def\<{\langle\,}
\def\>{\,\rangle}

\def\eg{{\em e.g. }}

\def\T{{\mathcal T}}
\def\F{{\mathcal F}}

\def\B{{\mathcal B}}

\def\g{\mathfrak g}
\def\nn{\mathfrak n}
\def\<{\langle}
\def\>{\rangle}
\def\n{{\mathfrak n}}

\def\1{\mathbf 1}
\def\Gr{{\rm Gr}}

\def\ii{{\mathbf i}}

\def\soc{{\rm Soc}}

\def\a{\alpha}

\def\L{\Lambda}
\def\mod{{\rm mod}}

\def\Ext{{\rm Ext}}
\def\add{{\rm add}\,}
\def\End{{\rm End}}
\def\Hom{{\rm Hom}}

\def\De{\Delta}

\def\la{\lambda}

\def\O{{\Omega}}
\def\1{{\mathbf 1}}
\def\ie{{\em i.e. }}

\def\G{\Gamma}

\def\b{\beta}
\def\ga{\gamma}

\def\bx{\mathbf x}

\def\bi{\mathbf{i}}

\def\Si{\Sigma}
\def\SS{\mathcal{S}}

\def\ds{\displaystyle}

\def\Ker{\mathrm{Ker\,}}

\def\bw{\overline{w}}
\def\bbw{\overline{\overline{w}}}
\def\bv{\overline{v}}
\def\bbv{\overline{\overline{v}}}
\def\Tr{\mathrm{T}}
\def\vpi{\varpi}

\def\stCC{\underline{\CC}}
\def\add{\mathrm{add}}
\def\G{\Gamma}

\def\RR{\mathcal{R}}
\def\ZZ{\mathcal{Z}}

\def\T{\mathcal{T}} 
\def\Sub{\mathrm{Sub}}
\def\md{\mathrm{mod}}
\def\dw{\dot{w}}
\def\tU{\widetilde{U}}
\def\tR{\widetilde{R}}
\def\l{\lambda}
\def\pr{{\rm pr}}
\def\vph{\varphi}

\begin{document}

\title{\bf Cluster structures on strata of flag varieties}

\author{B. Leclerc}
\date{}

\maketitle

\centerline{In memory of Andrei Zelevinsky, whose work has been a constant source of inspiration.}

\bigskip
\begin{abstract}
We introduce some new Frobenius subcategories of the module category of a preprojective algebra of Dynkin type, 
and we show that they have a cluster structure in the sense of Buan-Iyama-Reiten-Scott. 
These categorical cluster structures yield cluster algebra structures in the coordinate rings of  
intersections of opposite Schubert cells.
\end{abstract}

\section{Introduction}

\subsection{}
Let $G$ be a simple and simply connected algebraic group over $\C$. 
We assume that $G$ is simply-laced, that is, $G$ is of type $A, D, E$ in the Cartan-Killing classification.
We fix a maximal torus $H$ in $G$, a Borel subgroup $B$ containing $H$, and we denote by $B^-$ the Borel 
subgroup opposite to $B$ with respect to $H$. 
Let $W=\mathrm{Norm}_G(H)/H$ be the Weyl group, with length function $w\mapsto \ell(w)$
and longest element~$w_0$. 

We consider the flag variety $X = B^-\backslash G$, and we denote by $\pi:G\to X$ the natural projection
$\pi(g) := B^-g$. 
The Bruhat decomposition
\[
 G = \bigsqcup_{w\in W} B^-wB^-
\]
projects to the Schubert decomposition
\begin{equation}\label{eqSchubert}
X = \bigsqcup_{w\in W} C_w, 
\end{equation}
where $C_w = \pi(B^-wB^-)$ is the \emph{Schubert cell} attached to $w$, isomorphic to $\C^{\ell(w)}$.
We may also consider the Birkhoff decomposition
\[
 G = \bigsqcup_{w\in W} B^-wB,
\]
which projects to the \emph{opposite} Schubert decomposition
\begin{equation}\label{eqoppSchubert}
X = \bigsqcup_{w\in W} C^w, 
\end{equation}
where $C^w = \pi(B^-wB)$ is the \emph{opposite Schubert cell} attached to $w$, isomorphic to $\C^{\ell(w_0)-\ell(w)}$.
The intersection 
\begin{equation}
\RR_{v,w} := C^v\cap C_w
\end{equation}
has been considered by Kazhdan and Lusztig \cite{KL} in relation with the cohomological interpretation
of the Kazhdan-Lusztig polynomials.
One shows \cite{KL, D} that $\RR_{v,w}$ is non-empty if and only if $v\le w$ in the Bruhat order of $W$, 
and it is a smooth irreducible locally closed subset of $C_w$ of dimension 
$\ell(w)-\ell(v)$.
More recently, $\RR_{v,w}$ has
sometimes been called an open Richardson variety \cite{KLS}, because its closure in $X$ 
is known as a Richardson variety \cite{Rich}.

Intersecting the decompositions (\ref{eqSchubert}) and (\ref{eqoppSchubert}) of $X$, we thus get a finer stratification 
\begin{equation}\label{dec0}
X = \bigsqcup_{v\le w} \RR_{v,w}. 
\end{equation}
However, in contrast with (\ref{eqSchubert}) or (\ref{eqoppSchubert}), the strata $\RR_{v,w}$ of (\ref{dec0}) are not 
isomorphic to affine spaces. 

\subsection{}\label{sect-partial-flag}
Let $I$ be the vertex set of the Dynkin diagram of $G$. 
We denote by $x_i(t)\ (i\in I,\, t\in \C)$ (\resp $y_i(t)\ (i\in I,\, t\in\C)$) the one-parameter subgroups of $B$ (\resp $B^-$)
attached to the simple roots.
For $K\subset I$, let $B_K^-$ be the standard parabolic subgroup of $G$ generated by $B^-$
and the $x_i(t)$ with $i\in K$. We denote by $X_K=B_K^-\backslash G$ the corresponding
partial flag variety. Let $\pi_K: G \to X_K$ and $\pi^K: X\to X_K$ be the natural projections,
so that we have $\pi_K = \pi^K\circ\pi$.
Let $W_K$ be the parabolic subgroup of $W$ corresponding to $K$ with longest element $w_K$,
and let $W^K$ be the 
subset of $W_K$ consisting of the minimal length representatives of cosets in $W_K\backslash W$. 
It was shown by Lusztig \cite{L} that we have a decomposition
\begin{equation}\label{dec1}
 X_K = \bigsqcup_{w_Kv\le w,\ v\in W^K} \RR_{w_Kv,w}^K,
\end{equation}
where $\RR_{w_Kv,w}^K := \pi^K(\RR_{w_Kv,w})$ denotes the projection onto $X_K$ of the corresponding 
stratum of $X$. In fact $\pi^K$ restricts to an isomorphism 
$\pi^K\colon\RR_{w_Kv,w} \stackrel{\sim}{\to} \RR_{w_Kv,w}^K$, so the strata of (\ref{dec1}) are 
isomorphic to a subset of the strata of (\ref{dec0}).
Note also that (\ref{dec0}) can be regarded as a particular case of (\ref{dec1}) by choosing $K=\emptyset$.

\subsection{}
The stratification (\ref{dec1}) was introduced by Lusztig~\cite{L} for studying
the nonnegative part $(X_K)_{\ge 0}$ of the partial flag variety $X_K$. 
Lusztig conjectured that the intersections of the complex strata of (\ref{dec1}) with $(X_K)_{\ge 0}$ give
a cell decomposition of $(X_K)_{\ge 0}$, and this was proved by Rietsch \cite{R1}.
This cell decomposition was later studied in detail by  
Marsh-Rietsch \cite{MR}, Rietsch \cite{R}, and Rietsch-Williams \cite{RW}.
On the other hand, Goodearl and Yakimov \cite{GY1} have given a new interpretation of (\ref{dec1}) in a 
Poisson geometric setting. They showed that the strata of (\ref{dec1}) coincide with the $H$-orbits of the symplectic leaves
of the standard Poisson structure on $X_K$. 

The fact that the strata $\RR_{v,w}$ occur in these two settings naturally leads to the following question:
is there a \emph{cluster algebra structure} (in the sense of Fomin and Zelevinsky \cite{FZ2}) on the coordinate ring
$\C[\RR_{v,w}]$, which is compatible with total positivity and Poisson geometry?
More precisely, a cluster algebra structure on $\C[\RR_{v,w}]$ provides the complex variety
$\RR_{v,w}$ with a \emph{positive atlas}, and we would like each of the charts of this atlas to give us a description
of $\RR_{v,w} \cap (X)_{\ge 0}$. Moreover, if we denote by $\{\cdot,\cdot\}$ 
the Poisson bracket on $\C[\RR_{v,w}]$, we would like to have, for every cluster 
$\bx = \{x_1,\ldots,x_n\} \subset \C[\RR_{v,w}]$, the compatibility relations
\[
 \{x_i,x_j\} = \omega_{ij}x_ix_j,\qquad (1\le i<j \le n),
\]
for some appropriate $\omega_{ij}\in\Z$, see \cite[\S4.1]{GSV2}.

This situation has now become familiar because of several interesting classes of examples. 
The prototypical one, which has been one of the main motivations of Fomin and Zelevinsky for introducing 
cluster algebras, is given by the \emph{double Bruhat cells} in $G$ \cite{FZ1}. 
It was shown by Berenstein, Fomin and Zelevinsky \cite{BFZ} that the coordinate ring 
of a double Bruhat cell has the structure of an upper cluster algebra, which by construction
is compatible with the totally positive part described earlier in \cite{FZ1}. 
Recently, Goodearl and Yakimov \cite{GY2,GY3} have shown that the standard quantum deformation of the coordinate 
ring of a double Bruhat cell is a quantum cluster algebra confirming the Berenstein-Zelevinsky conjecture \cite{BZ2}. 
Passing to the classical limit, this shows that the coordinate ring is in fact a genuine cluster algebra.
Moreover, the cluster structure of \cite{BFZ} is compatible
with the standard Poisson structure on $G$ \cite{KZ}, \cite[\S4.3]{GSV2}. 
A second important class of examples, first studied by Gekhtman-Shapiro-Vainshtein \cite{GSV1}
and Fock-Goncharov \cite{FG}, is given by \emph{decorated Teichm\"uller spaces} of Riemann surfaces 
with marked points. Here the geometric object of interest, that is, the decorated Teichm\"uller space, appears as
the positive part of the cluster manifold of a cluster algebra defined in terms of triangulations
of the Riemann surface, and the compatible Poisson structure comes from the classical Weil-Petersson symplectic form.

\subsection{}\label{sect-catego-approach}
In this paper we will show the existence of a cluster structure 
for the strata $\RR_{v,w}$ of the flag variety $X$, and therefore also for the strata $\RR_{w_Kv,w}^K$ 
of the partial flag varieties~$X_K$.
We will follow the approach of the series of papers \cite{GLS2, GLS4, GLS5}, in which similar
questions were studied for an open cell of $X_K$, and for an arbitrary
unipotent cell of $G$ (or more generally, of a symmetric Kac-Moody group). 
The main feature of this approach is \emph{categorification}.
More precisely, let $\L$ be the preprojective algebra associated with $G$.
(For the definition of $\L$ and an introduction to its representation theory, see for instance
the survey papers \cite{GLS3,GLS6}.)
To each of the above varieties $Y$ (unipotent cell in $G$ or open cell in $X_K$), 
we have attached a certain Frobenius subcategory $\F_Y$ of $\mod(\L)$
having a cluster structure in the sense of Buan-Iyama-Reiten-Scott \cite{BIRS}.
This means that $\F_Y$ is endowed with a family of \emph{cluster-tilting objects}, which 
are related with each other by a procedure of \emph{mutation}.
Let $N$ be the unipotent radical of $B$. Using Lusztig's Lagrangian construction of the enveloping
algebra $U(\n)$ of $\n = \mbox{Lie}(N)$, we have constructed a \emph{cluster character} $\varphi$
from $\mod(\L)$ to $\C[N]$. The subspace spanned by the image of the restriction of $\varphi$ to $\F_Y$ is 
a certain subring of $\C[N]$, and we have shown that, after localizing with respect to the 
multiplicative subset given by the images of the projective objects of $\F_Y$, we obtain a ring
isomorphic to the coordinate ring of $Y$. Moreover, the images under $\varphi$ of the cluster-tilting
objects of $\F_Y$ endow this coordinate ring with the structure of a cluster algebra.

This categorical approach has several advantages. First, it gives a conceptual method to construct cluster 
algebra structures which might have been difficult to discover by elementary combinatorial methods,
see \eg \cite{GLS2}. Secondly, it provides a linear basis of the cluster algebra
coming from Lusztig's semicanonical basis of $U(\n)$. This basis contains all the cluster monomials.
Thirdly, the category $\F_Y$ has more structure than its geometric counterpart $Y$. In particular, 
it contains enough information to construct a quantum version of the cluster algebra structure on $\C[Y]$,
in the sense of \cite{BZ2}.
It was shown in \cite{GLS5} that this quantum cluster algebra is isomorphic to the quantized coordinate 
ring of $Y$ coming from the theory of quantum groups, and therefore by passing to the classical limit, the cluster variables
of any given cluster are compatible with the standard Poisson structure on~$\C[Y]$.

\subsection{}
We can now state the main result of this paper.
\begin{Thm}\label{mainThm1}
For every $v,w\in W$ with $v\le w$, the coordinate ring $\C[\RR_{v,w}]$  
contains a 
cluster subalgebra~$\tR_{v,w}$ coming from a Frobenius subcategory 
$\CC_{v,w}$ of the module category 
of the preprojective algebra $\L$. The number of cluster variables in a cluster of $\tR_{v,w}$
(including the frozen ones) is equal to $\dim(\RR_{v,w})$.
\end{Thm}

For a more detailed formulation of this result, see below Theorem~\ref{mainThm}.

If the category $\CC_{v,w}$ has finitely many indecomposable objects, then the cluster algebra $\tR_{v,w}$ has finite 
cluster type, and $\tR_{v,w} = \C[\RR_{v,w}]$, see \S\ref{sect-finite-type}. 
It is also easy to prove that $\tR_{v,w} = \C[\RR_{v,w}]$ if 
$w=v'v$ with $\ell(w) = \ell(v') + \ell(v)$, see \S\ref{sect-weak-Bruhat}.
We conjecture that this equality remains true in general.

Theorem~\ref{mainThm1} is new even in the most studied case of type $A_n$ Grassmannians $X_K = \Gr(j,n+1)$, that is, 
when $G=SL(n+1,\C)$, $K=[1,n]\setminus\{j\}$ and $v=w_Ku$ for some $u\in W^K$. 
In this case the strata 
$\RR_{w_Ku,w} \simeq\RR_{w_Ku,w}^K$, also called \emph{open positroid varieties}, have been studied by many authors, see 
\cite{P,BGY,KLS2, GLL1, GLL2}.
In their recent preprint \cite{MS}, Muller and Speyer have conjectured that the coordinate rings
of open positroid varieties are isomorphic to certain explicit cluster algebras described combinatorially
in terms of Postnikov diagrams. It is an interesting problem to compare this conjectural cluster structure
with the one given by Theorem~\ref{mainThm1} in this special case. A detailed example is studied in 
\S\ref{subsect-compare} below.

By construction $\CC_{v,w}$ is a subcategory of the Frobenius category
$\CC_w$ used in \cite{GLS4} to model the coordinate ring of the Schubert cell $C_w$. 
It then follows immediately from \cite{GLS5} that the cluster
algebra $\tR_{v,w}$ has a quantum deformation in the sense of \cite{BZ2}, which is a localization 
of a subalgebra of the quantum coordinate ring of the group $N(w)$ defined in Eq.~(\ref{defNw}) below. 
By passing to the classical limit, the cluster variables
of any given cluster of $\tR_{v,w}$ are compatible with the corresponding Poisson structure on $\C[\RR_{v,w}]$,
see \S\ref{sect-Poisson}.

We also believe that the positive atlas provided by $\tR_{v,w}$ gives the same positive
subset $\RR_{v,w}^{>0}$ as the one introduced by Lusztig, (see an example in \S\ref{example-positroid}). 
We hope to come back to this question in a forthcoming publication.

Here is a brief outline of the paper. 
In Section~\ref{sect3} we show that the coordinate ring of $\RR_{v,w}$ is isomorphic to a localization of an
invariant subring of $\C[N]$ (Theorem~\ref{thm-coordinate-ring}). 
This is needed to be in a situation where we can apply the categorification
approach of \S\ref{sect-catego-approach}, which only deals with (localizations of) subrings of~$\C[N]$.
In Section~\ref{sect2} we introduce the subcategory $\CC_{v,w}$ and show that it has a cluster structure in
the sense of \cite{BIRS}. We also give a general recipe for constructing cluster-tilting objects of~$\CC_{v,w}$, 
and we describe the $\CC_{v,w}$-projective-injective objects.
In Section \ref{sec-decategorif} we relate the category $\CC_{v,w}$ to the coordinate ring of 
$\RR_{v,w}$, and we prove our main result (Theorem~\ref{mainThm}). 
Section~\ref{sect-weak-Bruhat} deals with the special case $w=v'v$ with $\ell(w) = \ell(v') + \ell(v)$.
Section~\ref{sect-Poisson} discusses the quantization and the compatible Poisson structure of $\tR_{v,w}$.
Finally, \S\ref{sect5} illustrates the main ideas of the paper with examples.
In particular, for a specific open positroid variety in $\Gr(3,6)$, 
we compare in \S\ref{subsect-compare} the cluster algebra $\tR_{v,w}$ 
of Theorem~\ref{mainThm1} with the conjectural cluster structure of \cite[Conjecture 3.4]{MS}. 

\section{A quotient description of the strata}\label{sect3}

In this section we show that each stratum $\RR_{v,w}$ is isomorphic to the
quotient of an open subset $O_{v,w}$ 
of $N$ by the free action of the direct product $N(v)\times N'(w)$ of two subgroups of $N$.

\subsection{} 

For $w\in W$ define
\begin{equation}\label{defNw}
N(w) = N \cap (w^{-1}N^-w),\qquad  
N'(w) = N \cap (w^{-1}Nw),
\end{equation}
where $N^-$ denotes the unipotent radical of $B^-$.
As algebraic varieties, these unipotent subgroups are isomorphic to affine spaces: 
\[
N(w) \simeq \C^{\ell(w)}, \qquad N'(w) \simeq \C^{\ell(w_0)-\ell(w)}.
\]
We denote by $x\mapsto x^\Tr$ the ``transpose'' anti-automorphism of $G$ defined by 
\[
x_i(t)^\Tr = y_i(t),\quad
y_i(t)^\Tr = x_i(t),\quad
h^\Tr = h,\qquad
(i\in I,\ t\in\C,\ h\in H).
\]
We will make frequent use of the following well-known results, see \eg \cite[\S28]{H}.
\begin{Prop}\label{prop1} Let $w\in W$.
\begin{itemize}
 \item[(a)] The map $(x,y) \mapsto z=xy$ gives an isomorphism of algebraic varieties 
 $$N(w)\times N'(w) \stackrel{\sim}{\to} N.$$
 \item[(b)] Let $\dw$ be a fixed representative of $w$ in $\mathrm{Norm}_G(H)$.
 The map $(x,y) \mapsto z=x\dw y$ gives an isomorphism of algebraic varieties 
 $$B^-\times N'(w) \stackrel{\sim}{\to} B^-wB.$$
 \item[(c)] Similarly the map $(x,y) \mapsto z=x\dw y^{\Tr}$ gives an isomorphism of algebraic varieties 
 $$B^-\times N(w) \stackrel{\sim}{\to} B^-wB^-.$$
\end{itemize}
\end{Prop}
Inverting the isomorphism of Proposition~\ref{prop1} (b) and composing it with the projection on $N'(w)$
we get a well-defined morphism of algebraic varieties $z\mapsto y$ from $B^-wB$ to $N'(w)$
which we denote by $$z\mapsto [z]_w^+.$$
In other words, for $z\in G$, the intersection ${\dw}^{\,-1}B^-z \cap N'(w)$ is empty if $z\not \in B^-wB$, 
and it is reduced to a single point ${\dw}^{\,-1}B^-z \cap N'(w)= \{[z]_w^+\}$ if $z\in B^-wB$.
When $w=e$, the unit element of $W$, we have $N'(e) = N$ and we obtain a  
morphism $G_0 := B^-B \to N$ which we denote by $$z\mapsto [z]^+_e=[z]^+.$$

\subsection{}\label{sec-minors}
Following \cite[\S1.4]{FZ1}, to every $w\in W$ we attach two special representatives in $\mathrm{Norm}_G(H)$
denoted by $\bw$ and $\bbw$.
They satisfy \cite[Prop. 2.1]{FZ1}:
\begin{equation}
\bw^{\,-1} = \bw^{\Tr} = \overline{\overline{w^{-1}}}. 
\end{equation}
From now on, unless otherwise specified, the notation $\bw$ will be reserved for this special representative.
When an expression does not depend on the choice of a representative $\dw$ of $w$, we just write~$w$. 
For instance, the double coset $B^-\dw B^- \subset G$ does not depend on the choice of representative~$\dw$, 
so we just write $B^-wB^-$.

Let $\vpi_i\ (i\in I)$ be the fundamental weights of $\g={\rm Lie}(G)$. 
A regular function $\Delta^{\vpi_i}$ on $G$ is associated to each $\vpi_i$, as in \cite[\S1.4]{FZ1}.
For $u,v\in W$, one then defines following \cite[\S1.4]{FZ1} the generalized minor 
\begin{equation}\label{eqdefDelta}
 \Delta_{\,u(\vpi_i),\,v(\vpi_i)}(x) := \Delta^{\vpi_i}\left(\overline{\overline{u^{-1}}}x\,\overline{v}\right),\qquad (x\in G).
\end{equation}
This function only depends on the weights $u(\vpi_i)$ and $v(\vpi_i)$.

\subsection{} 
The following lemma shows that the strata $\RR_{v,w}$ are isomorphic to some
natural subvarieties of the unipotent subgroup $N$. 
It was proved by Brown, Goodearl and Yakimov \cite[Theorem 2.3(a)]{BGY}.
Since we use different notation and convention, we include a short proof
for the convenience of the reader.
For $v,w \in W$, define
\begin{equation}
 N_{v,w} := N'(v) \cap (v^{-1}B^-wB^-).
\end{equation}
\begin{Lem}[\cite{BGY}] \label{lemNwK}
The map $x\mapsto\pi(\bv x)$ induces an isomorphism of varieties $N_{v,w} \stackrel{\sim}{\longrightarrow} \RR_{v,w}$.
\end{Lem}

\proof
Let $v,w\in W$ with $v\le w$. We have by definition of $\RR_{v,w}$ 
\[
 \RR_{v,w} = C^v\cap C_w = \pi(B^-vB) \cap \pi(B^-wB^-) \supseteq \pi(B^-vB \cap B^-wB^-). 
\]
On the other hand, if $x\in \RR_{v,w}$, there exist $b_1,b_3,b_4\in B^-$ and $b_2\in B$
such that 
\[
x=\pi(b_1\overline{v}b_2) = \pi(b_3\bw b_4). 
\]
Therefore, by definition of $\pi$, there exists $b_5\in B^-$ such that
\[
b_1\overline{v}b_2 = b_5b_3\bw b_4 \in B^-vB \cap B^-wB^-, 
\]
so $x\in \pi(B^-vB \cap B^-wB^-)$. Hence $\RR_{v,w} = \pi(B^-vB \cap B^-wB^-)$.
Now by Proposition~\ref{prop1}(b), we have $B^-vB = B^-vN'(v)$ and $\pi$ restricts
to an isomorphism from $\bv N'(v)$ to $\pi(B^-vN'(v)) = C^v$.
It follows that $\pi$ restricts to an isomorphism from $\bv N'(v) \cap B^-wB^-$ to $\RR_{v,w}$.
Finally, left multiplication in $G$ by $\bv^{-1}$ restricts to an isomorphism from
$\bv N'(v) \cap B^-wB^-$ to $N'(v) \cap (v^{-1} B^-wB^-)$,
which finishes the proof.
\cqfd

\begin{remark}
{\rm
When $v=e$, we have $N'(v)=N$, and $N_{e,w} = N \cap B^-wB^-$
is a unipotent cell, studied in \cite{BZ1,GLS4}.
}
\end{remark}

\subsection{}
Let $O_{v,w}$ be the open subset of $N$ given by
\begin{equation}
 O_{v,w} := N \cap (v^{-1} G_0 w) = N \cap (v^{-1} B^-N w). 
\end{equation}
We have \cite[Prop. 7.3]{GLS4}:
\begin{equation}\label{eqdefOKw}
(x \in O_{v,w})  \Longleftrightarrow 
(x\in N \mbox{\ \ and\ \ } \bbv x\overline{w^{-1}}\in G_0)  \Longleftrightarrow  
(x\in N \mbox{\ \ and\ \ } \De_{v^{-1}(\vpi_i),w^{-1}(\vpi_i)}(x)\not = 0 \mbox{\ \ for all\ \ }i\in I).
\end{equation}
If $v \not \le w$ for the Bruhat order, then $v^{-1} \not \le w^{-1}$, and 
$v^{-1}(\varpi_i) \not \le w^{-1}(\varpi_i)$ for some $i\in I$
(for the induced Bruhat order on the set of extremal weights of the fundamental $\g$-module $L(\varpi_i)$,
isomorphic to $W/W_{I\setminus\{i\}}$).
Hence, for all $x\in N$ and some $i\in I$ we have 
$\De_{v^{-1}(\vpi_i),w^{-1}(\vpi_i)}(x) = 0$. Therefore if $v \not \le w$ then $O_{v,w} = \emptyset$.
From now on, unless otherwise specified, \emph{we will always assume that $v\le w$}.

\begin{Lem}\label{lemabc}
\begin{itemize}
\item[(a)] The subgroup $N(v)$ acts freely on $O_{v,w}$ by left multiplication.
\item[(b)] The subgroup $N'(w)$ acts freely on $O_{v,w}$ by right multiplication.
\item[(c)] We have $N_{v,w} \subset O_{v,w}$.
\end{itemize}
\end{Lem}

\proof
Let $n\in N(v)$. By definition of $N(v)$, we have $n\bv^{\,-1} = \bv^{\,-1}m$ for some $m\in N^-$. So
$nv^{-1} B^-N w \subseteq v^{-1}B^-Nw$. It follows that for every $x\in O_{v,w}$, we have $nx \in O_{v,w}$,
which proves (a).

\smallskip
Similarly, let $n'\in N'(w)$. By definition of $N'(w)$, we have $\bw n' = m'\bw$ for some $m'\in N$. So
$v^{-1} B^-N wn' \subseteq v^{-1}B^-Nw$. It follows that for every $x\in O_{v,w}$, we have $xn' \in O_{v,w}$,
which proves (b).

\smallskip
Let now $x\in N_{v,w}$. By Proposition~\ref{prop1}(c) we can write
$x = \bv^{\,-1} b \bw n^{\Tr}$ 
for some $b\in B^-$ and $n\in N(w)$.
Then $\bv x \bw^{\,-1} = b\bw n^{\Tr}\bw^{-1}$ and this is an element of $G_0$ because $n\in N(w)$.
This proves (c).
\cqfd

\begin{remark}
{\rm 
If $v\le w$, we know that $\RR_{v,w}\not = \emptyset$, hence, by Lemma~\ref{lemNwK}, 
$N_{v,w}\not = \emptyset$. 
So Lemma~\ref{lemabc}(c) implies that $O_{v,w}$ is a non-empty open subset of $N$, therefore
a dense open subset since $N$ is irreducible.
}
\end{remark}

\subsection{}
Define
\begin{equation}
 U_{v,w} := N \cap (B^-vBw^{-1}).
\end{equation}

\begin{Lem}\label{Lemme2.1}
The map $\zeta_{v,w} \ \colon\ x \mapsto [\bbv\,x\,\overline{w^{-1}}]^+$ is a surjective morphism from
$O_{v,w}$ to $U_{v,w}$. 
Moreover, $\zeta_{v,w}(x) = \zeta_{v,w}(x')$ if and only if there exists $n\in N(v)$ such that $x'=nx$.
\end{Lem}

\proof
Let $x\in N$. Then $x\in O_{v,w}$ if and only if $\bbv\,x\,\overline{w^{-1}} \in G_0$, hence $O_{v,w}$
is indeed the domain of $\zeta_{v,w}$. Moreover, 
\[
[\bbv\,x\,\overline{w^{-1}}]^+ = N\cap (B^-\bbv x \overline{w^{-1}})
\in N\cap (B^-v B w^{-1}) = U_{v,w}.
\]
Hence $\zeta_{v,w}(O_{v,w}) \subseteq U_{v,w}$. 
Conversely, if $y\in U_{v,w}$, there exist
$b\in B^-$ and $x\in N$ such that $y = b\bbv  x \overline{w^{-1}}$, and $y=\zeta_{v,w}(x)$.
Therefore $\zeta_{v,w} \colon O_{v,w} \to U_{v,w}$ is a surjective morphism.

If $x'=nx$ with $n\in N(v)$, we have $B^-\bbv x = B^-\bbv x'$. 
Indeed, 
\[
\bbv N(v) = \bbv N\cap N^-\bbv \subset B^-\bbv,
\]
hence $B^-\bbv x' = B^-\bbv nx \subseteq B^-\bbv x$, and similarly $B^-\bbv x' \subseteq B^-\bbv x$.
It follows that $\zeta_{v,w}(x) = \zeta_{v,w}(x')$.

Conversely, if $\zeta_{v,w}(x) = \zeta_{v,w}(x')$ there exist $n\in N$ and $b,b'\in B^-$ such that
\[
 n=b\bbv x\overline{w^{-1}} = b'\,\bbv x'\,\overline{w^{-1}},
\]
hence
\[
 x'x^{-1} = \bbv^{\,-1}(b')^{-1}b\bbv \in (v^{-1}B^-v)\cap N = N(v).
\]
\cqfd

\subsection{}
Define
\begin{equation}
 \O_{v,w} := N'(v) \cap O_{v,w}.
\end{equation}
By Proposition~\ref{prop1}(a), the multiplication map 
$N(v) \times N'(v) \stackrel{\sim}{\longrightarrow} N$ 
restricts to an isomorphism: 
\begin{equation} \label{eq11}
N(v) \times \O_{v,w} \stackrel{\sim}{\longrightarrow} O_{v,w}.
\end{equation}
Since $O_{v,w}$ is a dense open subset of $N$, it follows that $\O_{v,w}$ is a dense open subset of $N'(v)$.

\begin{Lem}\label{lemalg}
 The restriction of $\zeta_{v,w}$ to $\O_{v,w}$ is an isomorphism $\O_{v,w} \stackrel{\sim}{\longrightarrow} U_{v,w}$.
\end{Lem}

\proof
By Lemma~\ref{Lemme2.1} and (\ref{eq11}), the restriction of $\zeta_{v,w}$ to $\O_{v,w}$ is a bijective morphism
from $\O_{v,w}$ to $U_{v,w}$. It remains to show that its set-theoretic inverse is a morphism of algebraic
varieties.
Let $x\in \O_{v,w}$ and put $z=\zeta_{v,w}(x)$. By definition of $\zeta_{v,w}$,
there exists a unique $b\in B^-$ such that 
\[
\bbv x\overline{w^{-1}} = bz
\ \ \Longleftrightarrow \ \
z\bbw = b^{-1}\,\bbv x.
\]
Since $x\in\O_{v,w} \subset N'(v)$, it follows from Proposition~\ref{prop1}(b) that $x = [z\bbw]^+_v$.
Therefore the inverse map $z\mapsto x$ is obtained by composing the right multiplication by $\bbw$
with the morphism $y\mapsto [y]^+_v$, and 
$\zeta_{v,w}\colon \O_{v,w} \stackrel{\sim}{\longrightarrow} U_{v,w}$ is an isomorphism.
\cqfd

\subsection{}
Recall from Lemma~\ref{lemabc}(c) that $N_{v,w} \subset O_{v,w}$. Since by definition $N_{v,w} \subset N'(v)$
we have more precisely $N_{v,w} \subset \O_{v,w}$.

\begin{Lem}\label{lempenul}
The map $\zeta_{v,w}$ restricts to an isomorphism  
$\zeta_{v,w}\ \colon\  N_{v,w} \stackrel{\sim}{\longrightarrow} U_{v,w}\cap N(w^{-1})$.
\end{Lem}

\proof
Let $x\in N_{v,w}$. By Proposition~\ref{prop1}(c), we can write uniquely 
$x=\overline{v^{-1}}\, b\, \bbw \,n^{\Tr}$ with $b\in B^-$ and $n\in N(w)$.
It follows that 
$\zeta_{v,w}(x) = [\bbv x\overline{w^{-1}}]^+=[b\bbw n^{\Tr} \overline{w^{-1}}]^+$.
Since $n\in N(w)$, we have 
\[
n':=\bbw n^{\Tr} \overline{w^{-1}} \in N(w^{-1}), 
\]
hence
$\zeta_{v,w}(x) = n'$ belongs to $U_{v,w} \cap N(w^{-1})$.

Conversely, if $n'\in U_{v,w} \cap N(w^{-1})$, we have $n'\bbw = \bbw n^{\Tr}$ with $n\in N(w)$, hence
by the proof of Lemma~\ref{lemalg}, the unique preimage of $n'$ in $\O_{v,w}$ under $\zeta_{v,w}$ is 
$[n'\bbw]^+_v = v^{-1}B^-wn^{\Tr}\cap N'(v)\in N_{v,w}$. 
\cqfd

\begin{Lem}\label{lemfin}
\begin{itemize}
\item[(a)] The group $N'(w^{-1})$ acts freely by right multiplication on $U_{v,w}$. 
\item[(b)] Let $\phi_w(n) = \bbw n\overline{w^{-1}}$ be the conjugation group isomorphism from $N'(w)$ to $N'(w^{-1})$.
For $x\in O_{v,w}$ and $n\in N'(w)$ we have $\zeta_{v,w}(xn) = \zeta_{v,w}(x)\phi_w(n)$.
\end{itemize}
\end{Lem}

\proof
Let $x\in U_{v,w}$. This means that $x\in N$ and there exist $b^-\in B^-$ and $b\in B$
with $x=b^- \bv b \bw^{\,-1}$. Let $n\in N'(w^{-1})$. Then $\bw^{-1} n = n'\bw^{-1}$ with
$n'\in N'(w)$.
Hence $xn=b^- \bv bn' \bw^{\,-1} \in U_{v,w}$. This proves~(a). 
Then (b) is a straightforward calculation.
\cqfd

We are now ready for proving:
\begin{Thm}\label{theo2}
\begin{itemize}
 \item[(i)] The multiplication map $(n_1,y,n_2) \mapsto n_1yn_2$ gives an isomorphism $\mu_{v,w}$ from $N(v)\times N_{v,w}\times N'(w)$ to 
the open dense subset $O_{v,w}$ of $N$.
 \item[(ii)] Consider the invariant subring of $\C[O_{v,w}]$:
\[
 {}^{N(v)}\C[O_{v,w}]^{N'(w)}
:=\{f\in \C[O_{v,w}] \mid f(n_1yn_2) = f(y) \mbox{ for all } n_1\in N(v),\ n_2\in N'(w)\}.
\] 
The restriction of functions from $O_{v,w}$ to $N_{v,w}$ induces an algebra
isomorphism 
\[
{}^{N(v)}\C[O_{v,w}]^{N'(w)} \stackrel{\sim}{\longrightarrow} \C[N_{v,w}].
\]
\end{itemize}
\end{Thm}

\proof
The following commutative diagram, where the arrows labelled $\iota$ are inclusion maps, 
displays the morphisms appearing in the proof:
\[
 \xymatrix@-1.0pc{
&O_{v,w} \ar[rrd]^{\zeta_{v,w}}&&
\\
&\ar[u]^{\iota}\O_{v,w}\ar[rr]^{\zeta_{v,w}}&& U_{v,w}
\\
&\ar[u]^{\iota}N_{v,w}\ar[rr]^{\zeta_{v,w}} && U_{v,w} \cap N(w^{-1}) \ar[u]^\iota
}
\]
We first notice that, by Lemma~\ref{lemabc}, for $(n_1,y,n_2)\in N(v)\times N_{v,w}\times N'(w)$ we have 
$x=n_1yn_2 \in O_{v,w}$. Hence the image of $\mu_{v,w}$ is indeed contained in $O_{v,w}$.

Conversely, let $x\in O_{v,w}$.
By Eq.(\ref{eq11}) we can write $x = n_1\omega$, with $n_1\in N(v)$ and $\omega\in \O_{v,w}$.
By Lemma~\ref{lemalg} and Lemma~\ref{lemfin}(a) we can write $\zeta_{v,w}(\omega) = un'$,
with $u\in U_{v,w}\cap N(w^{-1})$ and $n'\in N'(w^{-1})$.
By the proof of Lemma~\ref{lemalg} we have $u=\zeta_{v,w}([u\bbw]^+_v)$, and $[u\bbw]^+_v \in N_{v,w}$
by Lemma~\ref{lempenul}. Hence
\[
 \zeta_{v,w}(\omega) = \zeta_{v,w}([u\bbw]^+_v) n' = \zeta_{v,w}([u\bbw]^+_v\phi_w^{-1}(n'))
\]
by Lemma~\ref{lemfin}(c). 
By Lemma~\ref{Lemme2.1} this implies that $\omega = n'_1[u\bbw]^+_v\phi_w^{-1}(n')$ 
for some $n'_1\in N(v)$. Hence, writing $n_2 = \phi_w^{-1}(n')$ and $y=[u\bbw]^+_v$, we have
\[
 x = n_1n'_1yn_2,
\]
where $n_1n'_1\in N(v)$, $y\in N_{v,w}$, and $n_2\in N'(w)$. Hence $\mu_{v,w}$ is surjective.

Now if
$x=n_1yn_2$ with $(n_1,y,n_2)\in N(v)\times N_{v,w} \times N'(w)$,
then by Lemma~\ref{Lemme2.1} and Lemma~\ref{lemfin}(c), we have 
\[
 \zeta_{v,w}(x)=\zeta_{v,w}(yn_2)= \zeta_{v,w}(y)\phi_w(n_2).
\]
Since $\zeta_{v,w}(y)\in N(w^{-1})$ and $\phi_w(n_2)\in N'(w^{-1})$, these two factors are uniquely
determined by $x$, and therefore $y$ and $n_2$ are uniquely determined by $x$. Hence $n_1$ is also 
uniquely determined by $x$. This proves that $\mu_{v,w}$ is injective, hence is a bijection 
from $N(v)\times N_{v,w}\times N'(w)$ to $O_{v,w}$.

Finally, $\mu_{v,w}^{-1}$ is a morphism of algebraic varieties. 
Indeed, by Proposition~\ref{prop1}(a), multiplication gives an isomorphism 
$N(w^{-1})\times N'(w^{-1}) \stackrel{\sim}{\to} N$. 
Let $\pr_1:N\to N(w^{-1})$ and $\pr_2:N\to N'(w^{-1})$ denote the two components of 
the inverse isomorphism. Let also $\xi_{v,w}$ denote the inverse of the isomorphism
of Lemma~\ref{lempenul}. Then if $x\in O_{v,w}$ and $(n_1,y,n_2) = \mu_{v,w}^{-1}(x)$  
we have
\[
 y=\xi_{v,w}\circ\pr_1\circ\zeta_{v,w}(x),\quad
 n_2=\phi_w^{-1}\circ\pr_2\circ\zeta_{v,w}(x),\quad
 n_1=xn_2^{-1}y^{-1},
\]
which shows that the maps $x\mapsto y$, $x\mapsto n_2$ and $x\mapsto n_1$ are 
morphisms.

Then (ii) follows immediately from (i).
\cqfd

\begin{remark}
{\rm Theorem~\ref{theo2}(i) can be rephrased by saying that
the action $(n_1,n_2)\cdot y = n_1yn_2^{-1}$ of $N(v)\times N'(w)$ on $O_{v,w}$ is free
with quotient $N_{v,w}$. However, the action of $N(v)\times N'(w)$ on $N$ is
\emph{not} free in general. 
In fact, $N(v)\times N'(w)$ acts freely on $N$ if and only if $w=v'v$ 
with $\ell(w)=\ell(v')+\ell(v)$, see below \S\ref{sect-weak-Bruhat}. 
}
\end{remark}

\subsection{}\label{sect-multiplicative}

By Eq.(\ref{eqdefOKw}), $O_{v,w}$ is the dense Zariski open subset of $N$ defined by the non-vanishing of 
\[
 D_{v,w}:=\prod_{i\in I} \De_{v^{-1}(\varpi_i),\,w^{-1}(\varpi_i)}. 
\]
Note that, as an element of $\C[N]$, the generalized minor $\De_{v^{-1}(\varpi_i),w^{-1}(\varpi_i)}$ is not always irreducible.
Let $M_{v,w}$ be the multiplicative subset of $\C[N]$ generated by the irreducible
factors of $D_{v,w}$.
It follows that the coordinate ring of
$O_{v,w}$ is obtained by localizing $\C[N]$ at $M_{v,w}$.
Therefore, using Theorem~\ref{theo2}, we see that $\C[N_{v,w}]$ is obtained by localizing the doubly-invariant
ring
\[
{}^{N(v)}\C[N]^{N'(w)}
:=\{f\in \C[N] \mid f(n_1xn_2) = f(x) \mbox{ for all } n_1\in N(v),\ n_2\in N'(w)\}
\]
at the same multiplicative subset $M_{v,w}$, which is clearly contained in ${}^{N(v)}\C[N]^{N'(w)}$.
In view of Lemma~\ref{lemNwK} we have thus proved:
\begin{Thm}\label{thm-coordinate-ring}
Let $S_{v,w}$ be the localization of ${}^{N(v)}\C[N]^{N'(w)}$ at $M_{v,w}$.
 The coordinate ring of the stratum $\RR_{v,w}$ 
 is isomorphic to $S_{v,w}$. \cqfd
\end{Thm}

\section{Categorical cluster structures}\label{sect2}

\subsection{}
We briefly recall some standard concepts of homological algebra. For more details see
\cite{ARS,ASS}. All algebras are finite-dimensional algebras over $\C$.

\subsubsection{}
Let $A$ be an algebra and
let $\CC$ be a subcategory of the category $\md(A)$ of finite-dimensional $A$-modules.
A \emph{left $\CC$-approxi\-ma\-tion} of $X\in\md(A)$ 
is a homomorphism $f\colon X \to Y$ with $Y\in \CC$, such that for every homomorphism $g\colon X\to Z$
with $Z\in \CC$ there exists $h\colon Y\to Z$ satisfying $g=hf$.
Dually, a \emph{right $\CC$-approximation} of $X$ 
is a homomorphism $f\colon Y \to X$ with $Y\in \CC$, such that for every homomorphism $g\colon Z\to X$
with $Z\in \CC$ there exists $h\colon Z\to Y$ satisfying $g=fh$.
We say that $\CC$ is \emph{functorially finite} if every $X\in\md(A)$ has a left $\CC$-approximation
and a right $\CC$-approximation.

Given $M\in\md(A)$, we denote by $\Sub(M)$ (\resp $\Fac(M)$) the full subcategory of $\md(A)$ whose objects
are the $A$-modules isomorphic to a submodule (\resp a factor module) of a direct sum of copies of $M$.
We denote by $\add(M)$ the additive closure of $M$, \ie the full subcategory of $\md(A)$ whose objects
are the $A$-modules isomorphic to a direct sum of direct summands of $M$.

We quote the following useful result of Auslander and Smal\o.

\begin{Prop}[{\cite[Theorem 5.10]{AS}}]\label{propAS}
Let $M$ and $M'$ be objects of $\md(A)$, and let $\tau$ be the Auslander-Reiten translation of $\md(A)$.
If 
\[
\Hom(\tau^{-1}(M), M) = 0 = \Hom(M',\tau(M')),
\]
then $\Sub(M) \cap\Fac(M')$ is an extension closed and functorially 
finite subcategory of $\md(A)$.
\end{Prop}

\subsubsection{}
Let $\CC$ be an extension closed full subcategory of $\md(A)$.
An $A$-module $M$ in $\CC$ is called $\CC$-{\it projective} (resp. $\CC$-{\it injective})
if $\Ext_A^1(M,X) = 0$ (resp. $\Ext_A^1(X,M) = 0$) for all $X \in \CC$.
If $M$ is $\CC$-projective and $\CC$-injective, then $M$ is also called $\CC$-{\it projective-injective}.
We say that $\CC$ has {\it enough projectives} (resp. {\it enough injectives})
if for each $X \in \CC$ there exists a short exact sequence
$0 \to Y \to M \to X \to 0$
(resp. $0 \to X \to M \to Y \to 0$)
where $M$ is $\CC$-projective (resp. $\CC$-injective) and $Y \in \CC$.
If $\CC$ has enough projectives and enough injectives, and if these coincide (i.e. an object is $\CC$-projective 
if and only if it is $\CC$-injective), then $\CC$ is called a {\it Frobenius category}.
If, moreover,  
for all $X,Y \in \CC$ there is a functorial isomorphism
\begin{equation}\label{functorial2}
\Ext_A^1(X,Y) \cong D\,\Ext_A^1(Y,X),
\end{equation}
where $D$ denotes the duality for vector spaces, we say that $\CC$ is \emph{stably $2$-Calabi-Yau}.
(The terminology comes from the fact that if $\CC$ is Frobenius, its stable category $\stCC$ is a triangulated
category; then (\ref{functorial2}) amounts to say that $\stCC$ is $2$-Calabi-Yau in the sense of Kontsevich.)

Frobenius stably 2-Calabi-Yau categories have been used by many authors to provide categorical models of cluster
algebras (see \eg \cite{GLS1, BIRS, FK}).

\subsection{}
Let $\L$ be the preprojective algebra over $\C$ of the same Dynkin type as $G$. 
We refer to \cite{GLS3,GLS6} for an introduction to this family of algebras and their representation theory. 
Here we will only recall their main features and review the results necessary 
for our current purpose.

\subsubsection{}
The $\C$-algebra $\L$ is finite-dimensional, basic, and selfinjective. 
Hence, $\md(\L)$ is an abelian Frobenius category. Moreover $\md(\L)$ is stably $2$-Calabi-Yau. 

The simple $\L$-modules are 1-dimen\-sio\-nal, in one-to-one correspondence with the vertices
of the Dynkin diagram of $G$.
We denote them by $S_i\ (i\in I)$, and we denote by $Q_i\ (i\in I)$ their injective envelopes.

Let $\tau$ be the Auslander-Reiten translation of $\md(\L)$. 
Because of the 2-Calabi-Yau symmetry, the $\tau$-translate of an indecomposable object $X\in\md(\L)$
is given by
\begin{equation}\label{tau}
\tau(X) = \Cok(X\to Q_X) 
\end{equation}
where $X\to Q_X$ is the injective envelope of $X$.
Dually, the $\tau^{-1}$-translate of $X$
is given by
\begin{equation}\label{tau-inv}
\tau^{-1}(X) = \Ker(P_X\to X) 
\end{equation} 
where $P_X\to X$ is the projective cover of $X$.

It was proved in \cite{GLS1} that $\md(\L)$ gives a categorical model for a cluster algebra structure
on the unipotent cell $N^{w_0} \subset N$, or equivalently on the open stratum $\RR_{e,w_0}$ 
of the flag variety $X$. 

\subsubsection{}\label{defEfunctors}

Denote by $s_i\ (i\in I)$ the Coxeter generators of $W$.
For $i\in I$, we define an endo-functor $\E_i=\E_{s_i}$ of $\md(\L)$ as follows \cite[\S5]{GLS2}.
Given $X\in\md(\L)$ we define $\E_i(X)$
as the kernel of a surjection 
$$
X \to S_i^{\oplus m_i(X)},
$$
where $m_i(X)$ denotes the multiplicity
of $S_i$ in the head of $X$.
If $f\colon X \to Y$ is a homomorphism,  
$f(\E_i(X))$ is contained in $\E_i(Y)$, and
we define 
$
\E_i(f)\colon \E_i(X) \to \E_i(Y)
$ 
as the restriction of $f$ to
$\E_i(X)$.
Clearly, $\E_i$ is an additive functor.
It acts on a module $X$ by removing the $S_i$-isotypical
part of its head.
Similarly, one defines a functor 
$\E^\dag_i=\E^\dag_{s_i}$ acting on $X$ by removing the 
$S_i$-isotypical part of its socle.

One shows that the functors $\E_i$ (\resp $\E_i^\dag$) satisfy the braid relations of $W$ \cite[Prop. 5.1]{GLS2}, 
hence by composing them we can define unambiguously functors $\E_w$ (\resp $\E_w^\dag$) for every $w\in W$. 

\begin{Lem}\label{lemE}
\begin{itemize}
 \item[(a)] If $v\le w$
in the Bruhat ordering, then $\E_w(X)$ is a submodule of $\E_v(X)$.
\item[(b)] If $X$ is a submodule of $Y$, then
$\E_w(Y/X) = \E_w(Y)/(X\cap \E_w(Y))$.
\end{itemize}
\end{Lem}
\proof
If $X$ is a submodule of $Y$ then we have a monomorphism $f: X \to Y$, which restricts to a monomorphism 
$\E_i(f): \E_i(X)\to \E_i(Y)$. Hence $\E_i(X)$ is a submodule of $\E_i(Y)$. This easily implies (a)
because we can write $w=s_{i_1}\cdots s_{i_r}$ and $v=s_{i_{j_1}}\cdots s_{i_{j_k}}$, for some
$1\le j_1<\cdots <j_k\le r$.

To prove (b), we first check it for $w=s_i$ and then use induction on $\ell(w)$. 
\cqfd

\subsubsection{}\label{defCw}
For $w\in W$, let $u := w^{-1}w_0$. Let $I_w := \E_u(\bigoplus_{i\in I} Q_i)$.
Following \cite{BIRS,GLS4}, we introduce the full subcategory
\begin{equation}
 \CC_w := \Fac(I_w)  
\end{equation}
of $\mod(\L)$.
One shows that $\CC_w$ is closed under extensions, Frobenius, and stably 2-Calabi-Yau, 
and that a module $X\in\CC_w$ is $\CC_w$-projective-injective if and only if it belongs to $\add(I_w)$
\cite{BIRS,GLS4}.
Moreover, we have 
\begin{equation}\label{eq-desc-Cw}
\CC_w = \E_{u}(\md(\L)).
\end{equation}
Indeed, if $Z = \E_u(Y)$ and $X$ is a submodule of $Z$, then $X$ is a submodule of $Y$ and we have
$Z/X = \E_u(Y/X)$ by Lemma~\ref{lemE}(b). 
It follows that $\E_{u}(\md(\L))$ is closed under taking quotient modules, and since
it contains $I_w$, it contains $\CC_w$. 
Conversely if $Z=\E_u(V)$ belongs to $\E_{u}(\md(\L))$, we can write $V=Y/X$ where $Y\in\add(\oplus_i Q_i)$
is the projective cover of $V$ in $\md(\L)$. Using Lemma~\ref{lemE}(b) again, we see that $Z$ is a quotient
of $\E_u(Y)\in \add(I_w)$, so $Z\in\CC_w$.

It was proved in \cite{GLS4} that $\CC_w$ gives a categorical model for a cluster algebra structure
on the unipotent cell $N_{e,w} \subset N$, or equivalently 
on the stratum $\RR_{e,w}$ of the flag variety $X$.

\subsubsection{}
Dually, for $w\in W$, let $J_w := \E^\dag_{w^{-1}}(\bigoplus_{i\in I} Q_i)$.
We introduce the full subcategory
\begin{equation}
 \CC^w := \Sub(J_w)  
\end{equation}
of $\mod(\L)$.
Again, $\CC^w$ is closed under extensions, Frobenius, stably 2-Calabi-Yau, 
its subcategory of $\CC^w$-projective-injective objects is equal to $\add(J_w)$, 
and we have 
\begin{equation}\label{eq-desc-Cv}
\CC^w = \E^\dag_{w^{-1}}(\md(\L)).
\end{equation}
Arguing as in \cite{GLS4}, it is not hard to show that $\CC^w$ gives a categorical model for a cluster algebra structure
on $N_{w,w_0} \subset N$, or equivalently 
on the stratum $\RR_{w,w_0}$ of the flag variety $X$. 
In the particular case when $w=w_K$ is the longest element of the parabolic subgroup $W_K$, the category $\CC^{w_K}$
was used in \cite{GLS2} to obtain a cluster structure on the open stratum $\RR_{w_K,w_0} \simeq \RR^K_{w_K,w_0}$ 
of the partial flag variety $X_K$
(see also \cite[\S17]{GLS4}, \cite{C1}).

\begin{example}
 {\rm
 Let $w=s_i$ be a simple reflexion. Then $I_w=S_i$ and $\CC_w$ is just the additive closure $\add(S_i)$ of $S_i$. 
 On the other hand, $\CC^{w} = \E_i^\dag(\md(\L))$ is the full subcategory 
 \[
  \CC^{s_i} = \{X\in\md(\L) \mid \Hom(S_i,X)=0\}.
 \]
In other words, $X\in\CC^{s_i}$ if and only if the socle of $X$ contains no copy of $S_i$.
 }
\end{example}

\subsubsection{}\label{sec-torsion-pair}
For every $w\in W$, the pair of subcategories $(\CC_w, \CC^w)$ is a \emph{torsion pair}.
This means \cite[VI.1]{ASS} that 
\begin{itemize}
 \item[(a)] $\Hom(M,N) = 0$ for every $M\in\CC_w$ and $N\in\CC^w$,
 \item[(b)] if $X\in\mod(\L)$ is such that $\Hom(X,N)=0$ for every $N\in\CC^w$, then $X\in\CC_w$,
 \item[(c)] if $Y\in\mod(\L)$ is such that $\Hom(M,Y)=0$ for every $M\in\CC_w$, then $Y\in\CC^w$.
\end{itemize}
This can be deduced for example from \cite[\S5]{BKT}, see in particular Example 5.14 and \S5.6.
(The categories $\CC_w$ and $\CC^w$ coincide with the categories denoted respectively by $\mathcal{T}_w$ 
and $\mathcal{F}^{w^{-1}w_0}$ in \cite{BKT}. Note that  
$\mathcal{F}^{w^{-1}w_0} = \mathcal{F}_w$ because $\L$ is of Dynkin type.) 

It follows that for every $X\in \md(\L)$, if we denote by $t_w(X)$ the maximal submodule of $X$ which is
contained in $\CC_w$, then we have $X/t_w(X) \in \CC^w$.
The assignment $X\mapsto t_w(X)$ defines a functor $t_w : \md(\L)\to \CC_w$ called the \emph{torsion radical} of 
the torsion pair $(\CC_w, \CC^w)$ \cite[VI.1.4]{ASS}.

Note that when $w=w_K=w_K^{-1}$ is the longest element of the parabolic subgroup $W_K$, we have
\cite[Proposition 5.4]{GLS2}
\begin{equation}
\E_{w_K}(X) = t_{w_0w_K}(X), \quad
\E^\dag_{w_K}(X) = X/t_{w_K}(X),\qquad (X\in\md(\L)).
\end{equation}
In general we only have a monomorphism and 
an epimorphism 
\begin{equation}\label{eqtw}
\E_{u}(X) \hookrightarrow  t_{w}(X),\quad
X/t_w(X) \twoheadrightarrow \E^\dag_{w^{-1}}(X),\qquad (X\in\md(\L)),
\end{equation}
where $u=w^{-1}w_0$.
Indeed, by (\ref{eq-desc-Cw}), we have that $\E_u(X)$ is a submodule of $X$ contained in $\CC_w$, 
so it is contained in $t_w(X)$.
Similarly, by (\ref{eq-desc-Cv}), we have that $\E^\dag_{w^{-1}}(X)$ is a quotient of $X$ contained in $\CC^w$, 
so it is a quotient of $X/t_w(X)$.
For instance, in type $A_2$ if $w=s_1s_2$ and $X=S_2$, we have $\E_{s_1s_2}(X)=0$
but $t_{s_2}(X) = X$.

However, when $X$ is a projective-injective object of $\md(\L)$, the morphisms (\ref{eqtw}) are isomorphisms,
as stated in the following
\begin{Lem}\label{lemtwproj}
Let $w\in W$ and set $u=w^{-1}w_0$. We have
\[
 t_w(Q_i) = \E_u(Q_i),\quad
 Q_i/t_w(Q_i) = \E^{\dag}_{w^{-1}}(Q_i), \qquad (i\in I).
\]
\end{Lem}
\proof
By (\ref{eqtw}) we have a short exact sequence
\[
 0\to \E_u(Q_i) \to t_w(Q_i) \to Y \to 0,
\]
with $Y= t_w(Q_i)/\E_u(Q_i)$. But $\E_u(Q_i)$ is $\CC_w$-projective-injective, and $Y\in\CC_w$ because
$\CC_w$ is closed under factor modules, so the sequence splits and $t_w(Q_i) \simeq \E_u(Q_i) \oplus Y$.
Since $Q_i$ has simple socle $S_i$, its submodule $t_w(Q_i)$ is indecomposable, so $Y=0$, and 
$t_w(Q_i) = \E_u(Q_i)$. The second equality is proved similarly.\cqfd

\begin{remark}\label{rem-inclusion-Cv}
{\rm
Note that if $v\le w$ then $v^{-1}w_0\ge w^{-1}w_0$, and by Lemma~\ref{lemE}(a), $I_v$ is a submodule
of~$I_w$. However, the fact that $v\le w$ does not imply that $\CC_v$ is contained in $\CC_w$. 
This is because a factor module of $I_v$ is not always a factor module of $I_w$.
This is only true if there exists $v'\in W$ such that
$w=v'v$ and $\ell(w) = \ell(v') + \ell(v)$. Indeed, in this case we have
\[
 \CC^w = \E^\dag_{w^{-1}}(\mod(\L)) = \E^\dag_{v^{-1}}(\CC^{v'}) \subseteq \E^\dag_{v^{-1}}(\md(\L)) = \CC^v, 
\]
hence $\CC^w\subseteq\CC^v$, and therefore $\CC_v \subseteq \CC_w$ because $(\CC_v,\CC^v)$ and $(\CC_w,\CC^w)$ are torsion pairs.
}
\end{remark}

\subsection{}
We retain the notation of the previous sections. Our main definition is
\begin{Def}
Let $v \le w$ in $W$.  We introduce the full subcategory
\[
 \CC_{v,w} := \CC^v \cap \CC_w
\]
of $\md(\L)$.
\end{Def}
Since $\CC^v$ and $\CC_w$ are closed under extensions, $\CC_{v,w}$ is closed under extensions.
Clearly, $\CC_{v,w}$ is also closed under direct sums and direct summands.
But it is in general neither closed under submodules nor under factor modules.

We now deduce from Proposition~\ref{propAS} that:
\begin{Prop}\label{propCKw}
The subcategory $\CC_{v,w}$ is functorially finite. 
\end{Prop}

\proof
Let $u = w^{-1}w_0$. We have $I_w = \bigoplus_{i\in I}\E_u(Q_i)$. 
Therefore, to show that $\Hom(I_w, \tau(I_w)) = 0$, it is enough to show that
\[
 \Hom(\E_u(Q_i), \tau(\E_u(Q_j)) = 0, \qquad (i,j\in I).
\]
If $\E_u(Q_j) \not = 0$, its injective envelope is $Q_j$, hence, by (\ref{tau}), $\tau(\E_u(Q_j))$ is isomorphic
to $Q_j/\E_u(Q_j)$. 
By Lemma~\ref{lemtwproj}, we have $Q_j/\E_u(Q_j)=Q_j/t_w(Q_j)\in\CC^w$, so indeed
$\Hom(\E_u(Q_i), \tau(\E_u(Q_j)) = 0$ because $(\CC_w,\CC^w)$ is a torsion pair.
Dually, a completely similar argument shows that
\[
\Hom(\tau^{-1}(J_v), J_v) = 0.
\]
The claimed result then follows from the definition
$\CC_{v,w} := \Sub(J_v) \cap \Fac(I_w)$
and Proposition~\ref{propAS}. \cqfd

\begin{remark}
{\rm
We can define $\CC_{v,w}$ even if $v\not\le w$. In fact, the proof of Proposition~\ref{propCKw} does 
not use the assumption $v\le w$. Similarly, Corollary~\ref{cor-clust}, Proposition~\ref{Prop-cluster-tilt}
and Proposition~\ref{prop-inj-proj-CKw} below are true for any pair $(v,w)\in W$. Moreover, 
$\CC_{v,w}=\{0\}$ if and only if $v = w'w$ with $\ell(v)=\ell(w')+\ell(w)$.
However, for our purpose
in this paper, we only need to consider the case $v\le w$.
}
\end{remark}

\subsection{}
We now review following \cite{BIRS} the notion of cluster structure of a subcategory of $\md(\L)$. 
Let $\B$ be a subcategory of $\md(\L)$ closed under extensions, direct sums and direct summands. 
For an object $T$ of $\B$
we denote by $\add(T)$ the \emph{additive envelope} of $T$, that is,
the full subcategory whose objects are finite direct sums
of direct summands of $T$.
We say that
$T$ is a $\B$-{\it cluster-tilting module} if 
\[
 (X \in \B \ \ \mbox{and} \ \ \Ext_\L^1(T,X) = 0) 
 \quad\Longleftrightarrow\quad
 (X \in \add(T)).
\]
We say that $T$ is \emph{basic} if its indecomposable direct summands are pairwise non-isomorphic.

Note that if $P$ is a $\B$-projective-injective module, then $P\in\add(T)$ for any $\B$-cluster-tilting module $T$. 
Note also that every $\B$-cluster-tilting module $T$ is \emph{rigid}, \ie $\Ext_\L^1(T,T)=0$.
In fact, since $\B$ is a subcategory of the module category of the preprojective algebra $\L$,
the module $T$ is $\B$-cluster-tilting if and only if it is a maximal rigid module in $\B$, \ie the number of
pairwise non-isomorphic indecomposable direct summands of $T$ is maximal among rigid modules in $\B$
\cite[Theorem II.3.1]{BIRS}.

\begin{Def}[\cite{BIRS}]\label{def-cluster-str}
The category $\B$ has a \emph{cluster structure} if
\begin{itemize}
 \item[(a)] $\B$ is a Frobenius category.
 \item[(b)] There exists a $\B$-cluster-tilting object.
 \item[(c)] If $T$ is a basic $\B$-cluster-tilting module and $M$ is an indecomposable direct summand of $T$ 
 which is not $\B$-projective-injective, then there exists a unique indecomposable object $M'\in\B$ such that 
 $T':=(T/M)\oplus M'$ is a $\B$-cluster-tilting module. Moreover we have nonsplit short exact sequences
 \[
   0\to M'\stackrel{f}{\longrightarrow} B \stackrel{g}{\longrightarrow} M\to 0,
   \qquad
   0\to M\stackrel{s}{\longrightarrow} B' \stackrel{t}{\longrightarrow} M' \to 0,
 \]
 where the maps $g$ and $t$ are minimal right $\add(T/M)$-approximations and $f$ and $s$ are minimal left
 $\add(T/M)$-approximations. We call $T'$ the mutation of $T$ at $M$, and we write $T'=\mu_M(T)$.
 \item[(d)] The Gabriel quiver $\G_T$ of the endomorphism algebra $\End(T)$ of any basic $\B$-cluster-tilting object $T$
 has no loop and no 2-cycle.
 \item[(e)] If $T'=\mu_M(T)$ then $\G_T$ and $\G_{T'}$ are related by a Fomin-Zelevinsky
 quiver mutation at the vertex of $\G_T$ corresponding to $M$. 
\end{itemize}
\end{Def}
We refer to \cite[\S3.2]{GLS1} for the definition of the quiver $\G_T$, and to \cite{FZ3} for the 
definition of quiver mutations.

It was shown in \cite{GLS1} that $\md(\L)$ has a cluster structure, and in \cite{BIRS,GLS4} that
$\CC_w$ (and dually $\CC^v$) has a cluster structure.
The following powerful criterion was proved by Buan, Iyama, Reiten and Scott.

\begin{Thm}[{\cite[Theorem II.3.1]{BIRS}}]\label{ThmBIRS}
Let $\B$ be an extension closed functorially finite subcategory of $\md(\L)$.
Then $\B$ has a cluster structure.
\end{Thm}

It follows from Proposition~\ref{propCKw} and Theorem~\ref{ThmBIRS} that:
\begin{Cor}\label{cor-clust}
 The subcategory $\CC_{v,w}$ has a cluster structure. \cqfd
\end{Cor}

\subsection{}\label{sect-recipe}
The following result gives a recipe  for constructing $\CC_{v,w}$-cluster-tilting objects.
\begin{Prop}\label{Prop-cluster-tilt}
Let $T$ be a $\CC_w$-cluster-tilting object. Then $T/t_v(T)$ is a $\CC_{v,w}$-cluster-tilting object. 
\end{Prop}
\proof
First note that $U:=T/t_v(T)\in\CC_{v,w}$. Indeed, $U\in\CC_w$ because $\CC_w$ is closed under factor modules,
and $U\in\CC^v$ because $t_v$ is the torsion radical of the torsion pair $(\CC_v,\CC^v)$.
Let $M\in\CC_{v,w}$. Then $M\in\CC_w$, and by \cite[Proposition 2.15]{GLS4} there exists an exact sequence
\begin{equation}\label{eqses}
 0\to T''\to T'\to M\to 0
\end{equation}
with $T'$ and $T''$ in $\add(T)$. 
Since $M\in \CC^v$, we have $t_v(M)=0$. It follows that $t_v(T')\simeq t_v(T'')$ and 
that we have a short exact sequence
\begin{equation}\label{eqses2}
 0\to T''/t_v(T'')\to T'/t_v(T')\to M\to 0.
\end{equation}
Write $U':= T'/t_v(T')$, and $U'':=T''/t_v(T'')$. 
Suppose that $\Ext_\L^1(U,M)=0$. Then, since $U''\in\add(U)$ the sequence (\ref{eqses2}) splits, 
hence $M$ is isomorphic to a direct summand of $U'$, so $M\in\add(U)$. 
Conversely, note that $T$ is rigid, that is, $T$ has an open orbit in its module variety \cite[Corollary 3.15]{GLS1}.
Using \cite[Lemma 4.3]{BKT}, this implies that $T/t_v(T)$ also has an open orbit in its module variety, that is, 
$U$ is also rigid. Hence if $M\in\add(U)$ then $\Ext_\L^1(U,M)=0$. 
Therefore $U$ is $\CC_{v,w}$-cluster-tilting.
\cqfd
\begin{remark}
{\rm
The module $T/t_v(T)$ is in general \emph{not} basic,
even if $T$ is a basic.
}
\end{remark}

We now describe the $\CC_{v,w}$-projective-injective objects.
\begin{Prop}\label{prop-inj-proj-CKw}
The module $M$ is $\CC_{v,w}$-projective-injective if and only if $M\in \add(I_w/t_v(I_w))$. 
\end{Prop}
\proof
We know that $I_w$ is a $\CC_w$-projective generator of $\CC_w$.
Therefore, if $M\in\CC_{v,w}\subset \CC_w$ we have an exact sequence
\[
 0\to M'' \to M' \to M \to 0
\]
where $M'\in \add(I_w)$ and $M''\in\CC_w$. Since $M\in\CC^v$ the same argument as in the proof of 
Proposition~\ref{Prop-cluster-tilt} shows that this yields an exact sequence
\[
 0 \to M''/t_v(M'') \to M'/t_v(M') \to M \to 0.
\]
If $M\not \in \add(I_w/t_v(I_w))$ then $M$ is not isomorphic to a direct summand of
$M'/t_v(M')$. So this last exact sequence does not split, hence $\Ext_\L^1(M,M''/t_v(M''))\not = 0$
and since $M''/t_v(M'')\in\CC_{v,w}$ this
shows that $M$ is not $\CC_{v,w}$-projective-injective.

Conversely, consider the short exact sequence
\[
 0 \to t_v(I_w) \to I_w \to I_w/t_v(I_w) \to 0.
\]
Let $M\in \CC_{v,w}$. Applying the functor $\Hom(-,M)$ to the above exact sequence, we get an exact sequence
\[
 0\to \Hom(I_w/t_v(I_w),M) \to \Hom(I_w,M) \to \Hom(t_v(I_w),M) \to \Ext^1(I_w/t_v(I_w),M) \to \Ext^1(I_w,M).
\]
But $\Ext^1(I_w,M) =0$, because $M\in \CC_w$ and $I_w$ is $\CC_w$-projective-injective.
Also, $\Hom(t_v(I_w),M) = 0$, because $t_v(I_w)\in \CC_v$ and $M\in \CC^v$. Therefore 
$\Ext^1(I_w/t_v(I_w),M) = 0$. This shows that $I_w/t_v(I_w)$ is $\CC_{v,w}$-projective-injective.
\cqfd

Recall that $I_w=\E_u(\oplus_{i\in I} Q_i)$, where $u=w^{-1}w_0$. 
Write $I_{i,w} := \E_u(Q_i)$, so that $I_w = \oplus_{i\in I} I_{i,w}$,
and 
\[
I_w/t_v(I_w) = \bigoplus_{i\in I} I_{i,w}/t_v(I_{i,w}).
\]
The next lemma gives a useful alternative description of the direct summands $Q_{i,v,w}:=I_{i,w}/t_v(I_{i,w})$.

\begin{Lem}\label{lem-proj-Cvw}
Let $v\le w$, and set $u=w^{-1}w_0$. We have
\[
Q_{i,v,w} =  \E^\dag_{v^{-1}}\E_u(Q_i),\qquad (i\in I).
\]
\end{Lem}
\proof
As was noted in Remark~\ref{rem-inclusion-Cv}, $I_{i,v}$ is a submodule of $I_{i,w}$.
Hence, using Lemma~\ref{lemtwproj}, 
\[
I_{i,v} = t_v(Q_i) = t_v(I_{i,w}).
\]
Now, again by Lemma~\ref{lemtwproj}, $Q_i/t_v(Q_i) = \E^\dag_{v^{-1}}(Q_i)$, so
$I_{i,w}/t_v(I_{i,w}) = \E^\dag_{v^{-1}}(I_{i,w}) = \E^\dag_{v^{-1}}\E_u(Q_i)$.
\cqfd

\begin{remark}\label{remark-proj-inj-CKw}
{\rm
The number of isomorphism classes of indecomposable $\CC_w$-projective-injective objects is at most $|I|$ \cite{BIRS,GLS4}.
But this is not always the case for $\CC_{v,w}$, because 
$Q_{i,v,w}$ is not always indecomposable, in contrast with $I_{i,w}$. 
This is illustrated in the following example. 

Let $\L$ be of type $A_3$, where we use the usual numbering of the Dynkin diagram by $I=\{1,2,3\}$. 
Let $w=s_1s_3s_2s_1s_3$ and $v=s_2$.
Then it is easy to check that there are 4 indecomposable $\CC_{v,w}$-projective-injective objects,
namely $Q_1, Q_3, S_1, S_3$. 
}
\end{remark}

\section{Cluster algebra structures}\label{sec-decategorif}

To obtain a cluster algebra from a cluster structure on a subcategory of $\md(\L)$ we use the 
cluster character $\vph : M\mapsto \varphi_M$ from $\md(\L)$ to $\C[N]$ described for example in \cite[\S6]{GLS4} 
(see also \cite{GLS3,GLS5}).  

\subsection{}
Let $M$ be a $\L$-module of dimension $d$. 
We quickly review the definition of $\varphi_M$. 
Given a finite sequence $\bi=(i_1,\ldots,i_e)\in I^e$, 
we introduce the variety $\F_{M,\bi}$ of
flags of submodules
\[
\mathfrak{f} = (0=F_0\subset F_1 \subset \cdots \subset F_e = M) 
\]
such that $F_k/F_{k-1} \cong S_{i_k}$ for $k=1,\ldots,e$.
This is a complex projective variety (empty if $e\not = d$).
We denote by $\chi(\F_{M,\bi})$ its topological Euler characteristic
(so in particular $\chi(\F_{M,\bi}) = 0$ if $e\not = d$).
The regular function $\varphi_M \in \C[N]$ is the unique polynomial function which
evaluates at an arbitrary product of one-parameter subgroups $x_j(t)$ of $N$ as
\begin{equation}
\vph_M(\mathbf{x}_{\mathbf{j}}(\mathbf{t})) 
= 
\sum_{\mathbf{a}\in\N^r} \chi(\F_{M,{\mathbf{j}}^{\mathbf{a}}}) \frac{\mathbf{t}^{\mathbf{a}}}{\mathbf{a}!},
\end{equation}
where $\mathbf{x}_{\mathbf{j}}(\mathbf{t}) := x_{j_1}(t_1)\cdots x_{j_r}(t_r)$, 
$\displaystyle\frac{\mathbf{t}^{\mathbf{a}}}{\mathbf{a}!} := \frac{t_1^{a_1}\cdots t_r^{a_r}}{a_1!\cdots a_r!}$,
and
${\mathbf{j}}^{\mathbf{a}} := (j_1,\ldots,j_1,j_2, \ldots,j_2, \ldots, j_r\ldots,j_r)$ with 
each component $j_k$ repeated $a_k$ times.
Thus the coefficients of $\varphi_M$ measure the ``size'' of the varieties of composition series of $M$,
where the ``measure'' of a variety is taken to be its Euler characteristic.
The existence of such a function $\vph_M$ follows from Lusztig's Lagrangian construction of $U(\n)$ 
\cite{L0,GLS1}.

\subsection{}
We now recall the main properties of the cluster character $\varphi$. 
They show that the ring $\C[N]$ can be regarded as a kind of dual Hall-Ringel algebra
(at $q=1$) of the category $\mod(\L)$. 

First, $\varphi$ maps direct sums to products, that is,
\begin{equation}\label{eq-cluster-char1}
 \varphi_{M\oplus M'} = \varphi_{M}\varphi_{M'},\qquad (M,M'\in\md(\L)).
\end{equation}
Next, if $\dim\Ext^1_\L(M,M')=1$, since $\md(\L)$ is stably 2-Calabi-Yau, we have nonsplit short exact sequences
\[
   0\to M'\to B \to M\to 0,
   \qquad
   0\to M \to B' \to M' \to 0,
\]
where $B$ and $B'$ are uniquely determined up to isomorphism. We then have \cite[Theorem 2]{GLSmult}
\begin{equation}\label{eq-cluster-char2}
 \varphi_{M}\varphi_{M'} = \varphi_{B} + \varphi_{B'}.
\end{equation}

It follows from these properties \cite{GLS1} that categorical mutations in the category
$\mod(\L)$ (in the sense of Definition~\ref{def-cluster-str}(c)) are transformed by $\vph$ into exchange relations 
in the Fomin-Zelevinsky cluster algebra structure of $\C[N]$, hence the name \emph{cluster character}.

\subsection{}
The generalized minors of \S\ref{sec-minors}, regarded as elements of $\C[N]$, are examples of functions $\varphi_M$.
Indeed we have (see \cite{GLS2}, \cite[Proposition 5.4]{C1}), if $v\le w$,
\begin{equation}\label{eq-phi-minors}
\varphi_{\E^\dag_v\E_{ww_0}(Q_i)} = \De_{v(\varpi_i),w(\varpi_i)},\qquad (i\in I). 
\end{equation}
Putting together Proposition~\ref{prop-inj-proj-CKw}, Lemma~\ref{lem-proj-Cvw} and (\ref{eq-phi-minors}), we then
obtain

\begin{Prop}\label{prop-phi-proj}
The category of $\CC_{v,w}$-projective-injective objects is equal to $\add(\oplus_{i\in I} Q_{i,v,w})$, and we
have 
\[
 \varphi_{Q_{i,v,w}} = \De_{v^{-1}(\varpi_i),w^{-1}(\varpi_i)},\qquad (i\in I).
\]
\cqfd
\end{Prop}

\subsection{}
Let $\B$ be a subcategory of $\md(\L)$ having a cluster structure.
We choose a basic $\B$-cluster-tilting object $T_0$ which we call \emph{initial}.
Another $\B$-cluster-tilting object $T$ is called \emph{reachable} if it can be obtained from $T_0$ via a finite 
sequence of mutations (it is then automatically basic). 
Let $\T$ denote the collection of all reachable $\B$-cluster-tilting objects.
Let $T=T_1\oplus\cdots\oplus T_k$ be a decomposition of $T\in\T$ into indecomposable direct summands.
We associate with $T$ the \emph{seed} 
\[
\Sigma_T:=(\G_T,(\varphi_{T_1},\ldots,\varphi_{T_k})).
\]
Let $R(\B)$ be the $\C$-subalgebra of $\C[N]$ generated by the functions $\varphi_{T_i}\ (1\le i\le k,\ T\in\T)$.  
It follows from Definition~\ref{def-cluster-str} and from (\ref{eq-cluster-char1}) (\ref{eq-cluster-char2})
that $R(\B)$ has the structure of a cluster algebra, with set of seeds $\{\Sigma_T\mid T\in\T\}$. 
The \emph{frozen} cluster variables (or \emph{coefficients}) of $R(\B)$ are the $\varphi_P$ where
$P$ is an indecomposable $\B$-projective-injective object.

Define also
 $S(\B) := \mbox{Span}_\C\langle \varphi_M \mid M\in \B\rangle.$
It follows from (\ref{eq-cluster-char1}) that $S(\B)$ is a subalgebra of $\C[N]$.
Clearly $R(\B) \subseteq S(\B)$.

\subsection{}
Let $v,w\in W$. It is shown in \cite{GLS4} that 
\begin{equation}\label{eqRCw}
R(\CC_w) = S(\CC_w)= \C[N]^{N'(w)} := \{f\in \C[N] \mid f(xn) = f(x) \mbox{ for all } x\in N,\ n\in N'(w)\}.
\end{equation}
Here, the initial cluster-tilting object used to define $R(\CC_w)$ is one of the explicit modules $V_{\bi}$
labelled by reduced words $\bi$ for $w$, whose definition is recalled in \S\ref{defVi} below. 
Dually, we have 
\begin{equation}\label{eqRCv}
R(\CC^v) =  S(\CC^v) = {}^{N(v)}\C[N] := \{f\in \C[N] \mid f(nx) = f(x) \mbox{ for all } x\in N,\ n\in N(v)\},
\end{equation}
(see \cite[\S17]{GLS4} and \cite{GLS2} for the particular case $v=w_K$ for some $K\subset I$).

\subsection{}
We recall the definition of the dual semicanonical basis of $\C[N]$. 
For more details, see \cite[\S2.2]{GLS4}.

Let $\bd\in\N I$ be a dimension vector for $\L$, and let $\L_{\bd}$ denote the affine variety of $\L$-modules
with dimension vector $\bd$. The cluster character $\varphi$ induces a constructible function from $\L_\bd$ to
$\C[N]$, that is, $\varphi_M$ takes only finitely many values on $\L_{\bd}$ and the fibres are constructible
subsets of $\L_{\bd}$. Hence if $Z$ is an irreducible component of the variety $\L_{\bd}$, there is a 
Zariski-dense subset of $Z$ over which $\varphi$ is constant. A module $M$ corresponding to a point
of this subset is called \emph{generic}, and we denote by $\varphi_Z$ the value $\varphi_M$ for a generic module
$M$ on $Z$. Let $\ZZ$ denote the union of all irreducible components of all varieties $\L_{\bd}\ (\bd\in \N I)$.
The collection $\SS^*:=\{\varphi_Z\mid Z\in\ZZ\}$ is dual to the semicanonical basis of $U(\nn)$ constructed
by Lusztig in \cite{L2}. Hence $\SS^*$ is a $\C$-basis of $\C[N]$ which we call the \emph{dual semicanonical basis},

Note that if $M$ is a rigid $\L$-module, then $M$ is generic and $\varphi_M \in \SS^*$.

\subsection{}
Recall from \S\ref{sect-multiplicative} the notation:
\begin{equation}
{}^{N(v)}\C[N]^{N'(w)} := \{f\in \C[N] \mid f(nxn') = f(x) \mbox{ for all } x\in N,\ n\in N(v), n'\in N'(w)\}.
\end{equation}

\begin{Prop}\label{prop-Svw}
We have:
\[
S(\CC_{v,w}) = {}^{N(v)}\C[N]^{N'(w)}.
\]
Moreover, the subset $\SS^* \cap S(\CC_{v,w})$ is a $\C$-basis of $S(\CC_{v,w})$.
\end{Prop}

\proof
Since $\CC_{v,w} = \CC^v\cap \CC_w$,
we have by (\ref{eqRCw}) and (\ref{eqRCv}),  
\[
S(\CC_{v,w}) \subseteq S(\CC^v) \cap S(\CC_w) = {}^{N(v)}\C[N] \cap  \C[N]^{N'(w)} = {}^{N(v)}\C[N]^{N'(w)}.
\]
Conversely, it is shown in \cite{GLS4} that $S(\CC_w)$ is spanned by a subset $\SS^*_w$ of the 
dual semicanonical basis $\SS^*$.
Moreover, if $Z\in\ZZ$ is the irreducible component corresponding to an element
$\varphi_Z$ of $\SS^*_w$, then there is an open dense subset of $Z$ consisting of modules $M\in\CC_w$,
see \cite[Proposition 14.5]{GLS4}.
Dually, $S(\CC^v)$ is spanned by a subset $\SS^{v,*}$ of $\SS^*$, and for every $\varphi_Z\in \SS^{v,*}$
there is an open dense subset of $Z$ consisting of modules $M\in\CC^v$.  

Hence the intersection $S(\CC^v) \cap S(\CC_w)$ is spanned by $\SS^{v,*} \cap \SS^*_w$, 
and if $\varphi_Z\in\SS^{v,*}\cap\SS^*_w$, then $\varphi_Z=\varphi_M$ for some 
$\L$-module $M\in \CC^v\cap\CC_w = \CC_{v,w}$. 
Therefore $S(\CC^v) \cap S(\CC_w)$ has a basis consisting of functions
$\varphi_M\in S(\CC_{v,w})$, so $S(\CC^v) \cap S(\CC_w) \subseteq S(\CC_{v,w})$.
\cqfd

\subsection{}
We now describe a cluster structure inside $S(\CC_{v,w})$.

\subsubsection{}\label{defVi}
Let $w=s_{i_r}\cdots s_{i_1}$ be a reduced decomposition, and put $\bi=(i_r,\ldots,i_1)$.
Let $R=\oplus_{i\in I}\Z\a_i$ denote the root lattice of $\g$, with basis the simple
roots $\a_i\ (i\in I)$. 
We identify dimension vectors $\bd = (d_i)$ for $\L$ with
elements of the positive cone $R^+=\oplus_{i\in I}\N\a_i$ via
\[
\bd \equiv \sum_{i\in I} d_i \a_i. 
\]
Consider the following sequence in $R^+$:
\[
\ga_k:=\varpi_{i_k} - s_{i_1}\cdots s_{i_{k}}(\varpi_{i_k}),\qquad  (1\le k \le r). 
\]
For $k = 1,\ldots, r$, one can show that there is a unique $V_k \in \md(\L)$ (up to isomorphism)
whose socle is $S_{i_k}$ and whose dimension vector is $\ga_k$.
Then \cite{GLS4} the module
\[
 V_\bi := \bigoplus_{k=1}^r V_k 
\]
is a cluster-tilting object of $\CC_w$.
Up to duality, this $\L$-module is the same as the one introduced
by Buan, Iyama, Reiten, and Scott \cite{BIRS}.

\subsubsection{}\label{sect-finite-type}
We can now use the recipe of \S\ref{sect-recipe}.
Define $U_\bi := V_\bi/t_v(V_\bi)$. By Proposition~\ref{Prop-cluster-tilt}, $U_{\bi}$ is a cluster-tilting
object in $\CC_{v,w}$.
Using $U_\bi$ as an initial cluster-tilting object in $\CC_{v,w}$, we can then define the cluster algebra 
\[
R(\CC_{v,w}) \subseteq S(\CC_{v,w}).
\]
All cluster monomials of $R(\CC_{v,w})$ belong to the dual semicanonical basis $\SS^*\cap S(\CC_{v,w})$ 
of $S(\CC_{v,w})$.

If the category $\CC_{v,w}$ has a finite number of isomorphism classes of indecomposable objects, that is,
if the cluster algebra $R(\CC_{v,w})$ has finite cluster type (see \cite[\S18]{GLS3}), then every indecomposable
module $X\in\CC_{v,w}$ is rigid and gives rise to a cluster variable $\varphi_X \in R(\CC_{v,w})$. It follows
from (\ref{eq-cluster-char1}) that for every $M\in\CC_{v,w}$ the regular function $\varphi_M$ is a product of cluster 
variables, and therefore $R(\CC_{v,w}) = S(\CC_{v,w})$.

\subsubsection{}\label{defvj}

We now give a more explicit description of the cluster-tilting module $U_\bi$.
Inspired by 
\cite[\S3]{MR}, we define inductively a sequence $(v_{(0)}, v_{(1)}, \ldots, v_{(r)})$
of elements of~$W$, by setting $v_{(0)}=e$, and for $k=1,\ldots,r$,
\[
 v_{(k)} = 
 \left\{
 \begin{matrix}
s_{i_k}v_{(k-1)}&\mbox{if\ \ } vv_{(k-1)}^{-1}s_{i_k} < vv_{(k-1)}^{-1},\\[2mm]
v_{(k-1)}&\mbox{otherwise.}
 \end{matrix}
\right.
 \]
It is easy to check that $v_{(0)}\le v_{(1)}\le \cdots \le v_{(r)} = v$.
Let 
\[
 J_{v,\bi} := \{k\in\{1,\ldots, r\} \mid v_{(k)} = v_{(k-1)}\}.
\]
Since $r=\ell(w)$ and $\ell(v_{(k)})=\ell(v_{(k-1)})+1$ if $k\not\in  J_{v,\bi}$, we see that
\[
| J_{v,\bi}| = \ell(w)-\ell(v) = \dim \RR_{v,w}.
\]
Define the $\L$-modules
\begin{equation}
 U_j := \E^\dag_{v_{(j)}^{-1}} V_j,\qquad (1\le j\le r).
\end{equation}
Finally, for a $\L$-module $M$, denote by $\Si(M)$ the number of pairwise non-isomorphic indecomposable 
direct summands of $M$.
\begin{Prop}\label{prop-desc-U}
We have $U_j = V_j/t_v(V_j)$, and therefore $U_\bi = \bigoplus_{1\le j\le r} U_j$. 
Moreover, 
\[
\Si(U_\bi) = \ell(w)-\ell(v).
\]
\end{Prop}

\proof
Let us recall from \cite[\S2.4]{GLS4} the description of the direct summands $V_s\ (1\le s\le r)$ of 
the $\CC_w$-cluster-tilting module $V_\bi$.
For a $\L$-module $X$ and a simple module $S_j$, let
$\soc_{(j)}(X)$ be the sum of all submodules
$U$ of $X$ with $U \cong S_j$.
For a sequence $(j_1,\ldots,j_t)\in I^t$, there is a unique sequence
$$
0 = X_0 \subseteq X_1 \subseteq \cdots \subseteq X_t \subseteq X
$$
of submodules of $X$ such that $X_p/X_{p-1} = \soc_{(j_p)}(X/X_{p-1})$.
Define $\soc_{(j_1,\ldots,j_t)}(X) := X_t$.
Then, if  $\bi = (i_r,\ldots,i_1)$ is a reduced expression of $w$,
we have 
$$
V_s = \soc_{(i_s,\ldots,i_1)}\left(Q_{i_s}\right),\qquad (1\le s\le r), 
$$
and $V_\ii = V_1 \oplus \cdots \oplus V_r$.

Let $K_{v,\bi}=\{k_1<k_2<\cdots<k_l\}$ be the complement of $J_{v,\bi}$ in $\{1,\ldots,r\}$.
By construction of the sequence $(v_{(k)})$, we have $v=s_{i_{k_l}}\cdots\, s_{i_{k_2}} s_{i_{k_1}}$,
and $l=\ell(v)$. 
Given $s\in\{1,\ldots,r\}$ we have 
$v_{(s)} = s_{i_{k_t}}\cdots\, s_{i_{k_2}} s_{i_{k_1}}$ for some $t\le s$.
Define
\[
W_s = \soc_{(i_{k_t},\ldots,i_{k_1})}\left(Q_{i_{s}}\right),\qquad (1\le s\le r).
\]
Then each module $W_s$ is a direct summand of a $\CC_v$-cluster-tilting module, and in particular
belongs to $\CC_v$. Moreover, again by construction of the sequence $(v_{(k)})$, 
$W_s$ is the largest submodule of $V_s$ contained in $\CC_v$. Hence
\[
 U_s = \E^\dag_{v_{(s)}^{-1}} V_s = V_s/W_s = V_s/t_v(V_s).
\]
This proves the first statement. 

For the second one, we recall from \cite[\S10]{GLS4} that each module $Y\in\CC_w$ has a canonical
filtration 
\[
0 = Y_0 \subseteq Y_1 \subseteq \cdots \subseteq Y_r \subseteq X 
\]
such that $Y_k/Y_{k-1} \simeq M_k^{a_k}$ for some $a_k\in \N$.
Here $M_k = V_k/V_{k^-}$, where $k^-=\max\{s\le k-1 \mid i_s = i_k\}$.
Now one can check that, if $s\in J_{v,\bi}$, 
the upper layer of the canonical filtration of $Y=U_s$
is equal to $M_s$, that is, $U_s$ has the same top layer as the indecomposable module $V_k$. 
This proves that the module $\bigoplus_{s\in J_{v,\bi}} U_s$ contains
at least $|J_{v,\bi}|=\ell(w)-\ell(v) = \dim \RR_{v,w}$
pairwise non-isomorphic indecomposable direct summands. But because of Theorem~\ref{theo2} and Proposition~\ref{prop-Svw}, 
a basic rigid object 
\[
T= T_1\oplus\cdots\oplus T_m \in \CC_{v,w}
\]
has at most $\dim \RR_{v,w}$ 
pairwise non-isomorphic direct summands. Indeed, all the monomials $\varphi_{T_1}^{b_1}\cdots \varphi_{T_m}^{b_m}$
belong to the dual semicanonical basis $\SS^*$, hence are linearly independent. So the regular functions
$\varphi_{T_1},\ldots,\varphi_{T_m}$ are algebraically independent in $\C[\RR_{v,w}]$.
This forces 
\[
\Si(U_\bi) = \Si\left(\bigoplus_{s\in J_{v,\bi}} U_s\right) = \ell(w)-\ell(v).
\]
\cqfd

\subsubsection{}
We can now describe the initial cluster of the cluster algebra $R(\CC_{v,w})$ given by the cluster-tilting 
object $U_\bi$.

\begin{Cor}\label{cor-number}
For $j\in\{1,\ldots,r\}$, set $w_{(j)}=s_{i_j}\cdots s_{i_1}$.
Recall also from \S\ref{defvj} the elements $v_{(j)}\in W$.
The irreducible factors of the generalized minors
\[
 \De_{v_{(j)}^{-1}(\varpi_{i_j}), w_{(j)}^{-1}(\varpi_{i_j})},\qquad (j\in J_{v,\bi}), 
\]
are the $\ell(w)-\ell(v)$ cluster variables of the initial cluster of $R_{v,w}$ corresponding 
to $U_\bi$. The frozen cluster variables are those which divide
$\prod_{i\in I} \De_{v^{-1}(\varpi_{i}), w^{-1}(\varpi_{i})}$.
\end{Cor}

\proof
This follows immediately from (\ref{eq-phi-minors}), Proposition~\ref{prop-desc-U} and 
Proposition~\ref{prop-phi-proj}.
\cqfd

\subsubsection{}
By Proposition~\ref{prop-phi-proj}, the multiplicative subset of $\C[N]$ generated by the regular functions $\varphi_P$ where
$P$ runs over the set of $\CC_{v,w}$-projective-injective objects is nothing else than the multiplicative subset $M_{v,w}$
defined in \S\ref{sect-multiplicative}. Therefore, putting together Theorem~\ref{thm-coordinate-ring},  
Proposition~\ref{prop-Svw}, and Corollary~\ref{cor-number}, we see that we have now proved Theorem~\ref{mainThm1}. 
Here is a more precise formulation of the results we have obtained.

\begin{Thm}\label{mainThm}
\begin{itemize}
\item[(i)]For every $v,w\in W$ with $v\le w$, there is a Frobenius subcategory $\CC_{v,w}$ of 
$\md(\L)$ such that the subspace $S(\CC_{v,w})$ of $\C[N]$ spanned by the 
functions $\varphi_M\ (M\in \CC_{v,w})$ is equal to the doubly invariant ring ${}^{N(v)}\C[N]^{N'(w)}$.
\item[(ii)]
The localization $S_{v,w}$ of $S(\CC_{v,w})$
with respect to the multiplicative subset 
\[
M_{v,w}=\{\varphi_P \mid P \mbox{ projective in } \CC_{v,w}\} 
\]
is isomorphic to $\C[\RR_{v,w}]$.
\item[(iii)] The category $\CC_{v,w}$ has a cluster structure in the sense of {\rm\cite{BIRS}}.
The subalgebra $R(\CC_{v,w})$ of $S(\CC_{v,w})$ generated by the $\varphi_M$ where $M$ runs
over the set of direct summands of reachable cluster-tilting objects of $\CC_{v,w}$, has the 
structure of a cluster algebra.
The set $M_{v,w}$ coincides with the set of monomials in the frozen variables of $R(\CC_{v,w})$.
Hence, the cluster algebra $R_{v,w}$ obtained from $R(\CC_{v,w})$ by inverting the frozen cluster
variables is isomorphic to a cluster subalgebra $\tR_{v,w}$ of $\C[\RR_{v,w}]$.

\item[(iv)] The number of cluster variables in a cluster of $R_{v,w}$ (including the frozen ones) is
equal to $\dim\RR_{v,w}$.

\item[(v)]  By \S\ref{sect-finite-type}, when $\CC_{v,w}$ has finitely many indecomposable objects up to isomorphism, we have
$R(\CC_{v,w}) = S(\CC_{v,w})$ and the cluster algebra $\tR_{v,w}$ is equal to  $\C[\RR_{v,w}]$.
\end{itemize}
\end{Thm}

Note that we have also proved (Proposition~\ref{prop-Svw}) that $S(\CC_{v,w})$ is spanned by a subset of the 
dual semicanonical basis $\SS^*$. Multiplying the elements of this subset by the inverses of the elements
of $M_{v,w}$, we thus obtain a $\C$-basis of $S_{v,w}\simeq\C[\RR_{v,w}]$ which contains all the cluster monomials
of $R_{v,w}$.

In type $A_n$ with $n\le 4$, the category $\md(\L)$ has finitely many indecomposables, so
Theorem~\ref{mainThm}~(v) applies for every pair $(v,w)$ in this case.
More generally, consider a parabolic subgroup $W_K$ of $W$, and let $w_K$ be its longest element.
The categories $\CC^{w_K}$ with finitely many indecomposables have been classified in \cite[\S11]{GLS2}.
We have $\CC^{w_Kv} = \E^\dag_{v^{-1}}(\CC^{w_K})$, so if $\CC^{w_K}$ belongs to the
list of \cite[\S11]{GLS2}, then $\CC_{w_Kv,w}$ has finitely many indecomposables for every
pair $v\in W^K$ and every $w\ge w_Kv$. Thus, by Eq.(\ref{dec1}), the coordinate rings of all Lusztig strata of the partial
flag variety $X_K$ are cluster algebras in this case.

\section{The case $w=v'v$ with $\ell(w)=\ell(v')+\ell(v)$}\label{sect-weak-Bruhat}

In this section we consider the property: 
\begin{equation*}\label{property}
\mbox{(P)\quad $w\in W$ can be factored as $w=v'v$ with $\ell(w)=\ell(v')+\ell(v)$.} 
\end{equation*}

\subsection{}
Property (P) implies that $v\le w$ in the Bruhat ordering, but, as is well known, it is 
stronger.
In fact, (P) holds if and only if $N(v)\subseteq N(w)$, see \eg \cite[\S2]{BZ1}.
In other words, (P) is equivalent to 
\[
N(v)\cap N'(w) = \{e\}.
\]
This is in turn equivalent to the fact that $N(v)\times N'(w)$ acts freely on the whole group $N$.
Indeed, if $n_1xn_2 = x$ with $n_1\in N(v)$, $x\in N$, $n_2\in N'(w)$, then 
writing $n=n'n''$ with $n'\in N(v)$ and $n''\in N'(v)$ we obtain
\[
 (n_1n')(n''n_2) = n'n''.
\]
Here $n_1n'\in N(v)$ and $n''n_2\in N'(v)$ since $N'(w)\subseteq N'(v)$, so
this forces $n_1n'=n'$ and $n''n_2=n''$, hence $n_1=n_2=e$.
Conversely, since $N(v)\cap N'(w)$ is isomorphic to the stabilizer of $e$, if (P) does not hold the action is not free.  

So let us assume that (P) holds. A natural fundamental domain for the action 
of $N(v)\times N'(w)$ on $N$ is the unipotent subgroup $N_v^w$ corresponding
to the subset of positive roots
\[
 \De_v^w = \{\b\in \De^+ \mid w(\b)<0,\ v(\b)>0\}. 
\]
Thus, $N(v)\backslash N/N'(w) \simeq N_v^w$ is an affine space of dimension $\ell(w)-\ell(v)$.

\subsection{}
By Remark~\ref{rem-inclusion-Cv}, property (P) implies (and actually is equivalent to) 
$\CC_v \subseteq \CC_w$ and $\CC^w \subseteq \CC^v$. 
In this situation, we write $(\CC_v,\CC^v)\preccurlyeq (\CC_w,\CC^w)$.
Such nested torsion pairs have been studied extensively in \cite{BKT} in relation with
Harder-Narasimhan filtrations and stability conditions.

Write $q=\ell(v)$ and $r=\ell(w)$. 
If (P) holds, we can choose a reduced word $\bi = (i_r,\ldots,i_1)$ for $w$ such that $(i_q,\ldots,i_1)$
is a reduced word for $v$. Clearly we then have $J_{v,\bi} = \{q+1,q+2,\ldots,r\}$, and
\[
\tU_\bi:= U_{q+1}\oplus \cdots \oplus U_r = V_{q+1}/t_v(V_{q+1}) \oplus \cdots \oplus  V_{r}/t_v(V_{r}) 
= \E^\dag_{v^{-1}}(V_{q+1}) \oplus \cdots \oplus \E^\dag_{v^{-1}}(V_{r})
\]
is a basic $\CC_{v,w}$-cluster-tilting module.

Recall from the proof of Proposition~\ref{prop-desc-U} the modules $M_k\ (1\le k\le r)$ in $\CC_w$. 
By \cite[Lemma 9.8(i)]{GLS4}, we have $\Hom(M_k,M_l)=0$ if $k<l$.
It then follows from \cite[Theorem 5.11(iii)]{BKT} that $M_l\in\CC^v$ for every $l = q+1,\ldots,r$.
Therefore property (P) implies that 
\[
 M_l \in \CC_{v,w},\qquad (q+1\le l \le r). 
\]
Thus, using \cite[Theorem 15.1]{GLS4}, the monomials in the functions $\varphi_{M_l}\ (q+1\le l \le r)$
form a dual PBW-basis of $\C[N_v^w]$. Note in particular that the set of dimension vectors
of the modules $M_l\ (q+1\le l \le r)$ coincides with $\De_v^w$.
Moreover a straightforward adaptation 
of \cite[\S13]{GLS4} gives a mutation sequence in $\CC_{v,w}$ from $\tU_\bi$ to 
$\O_{v,w}^{-1}(\tU_\bi)$ such that all modules $M_l\ (q+1\le l \le r)$ occur once in this sequence.
(Here $\O_{v,w}$ denotes the syzygy functor of the Frobenius category $\CC_{v,w}$.)
Therefore all the functions $\varphi_{M_l}\ (q+1\le l \le r)$ are cluster variables of
$R(\CC_{v,w})$, and this implies that $R(\CC_{v,w})=\C[N_v^w]$.
Finally, localizing with respect to the multiplicative subset $M_{v,w}$, we have shown:
\begin{Prop}
If property (P) holds, then the cluster algebra $\tR_{v,w}$ is equal to  
$\C[\RR_{v,w}]$.\cqfd
\end{Prop}
\begin{remark}
{\rm
By \cite[Proposition 5.16]{BKT}, if Property (P) holds we have an equivalence of categories
$\CC_{v'}\stackrel{\sim}{\to}\CC_{v,w}$ via a certain reflection functor $\Si$.
In fact the $\CC_{v,w}$-cluster-tilting module $\tU_\bi$ is the image by $\Si$ of the
$\CC_{v'}$-cluster-tilting module associated with the reduced word $(i_r,\ldots,i_{q+1})$.
}
\end{remark}

\section{Quantization and Poisson structure}\label{sect-Poisson}

\subsection{} Let $w\in W$.
In \cite[\S10.7]{GLS5}, the cluster algebra structure on $\C[N(w)] \simeq \C[N]^{N'(w)}$ 
given by the subcategory~$\CC_w$ has been quantized into a $\Q[q^{\pm1}]$-algebra $\A_{\Q[q^{\pm1}]}(\CC_w)$.
For every reachable basic cluster-tilting object $T=T_1\oplus\cdots \oplus T_r$ of $\CC_w$,
the associated cluster $(\varphi_{T_1},\ldots,\varphi_{T_r})$ in $\C[N]^{N'(w)}$ is quantized into a
collection $(Y_{T_1},\ldots,Y_{T_r})$ of elements of $\A_{\Q[q^{\pm1}]}(\CC_w)$ satisfying the $q$-commutation relations
\begin{equation}\label{eq-q-commut}
 Y_{T_i}Y_{T_j} = q^{\l_{ij}} Y_{T_j}Y_{T_i},\qquad (1\le i,j \le r),
\end{equation}
where
\begin{equation}\label{log-can}
 \l_{ij} = \dim\,\Hom_\L(T_i,T_j) - \dim\,\Hom_\L(T_j,T_i).
\end{equation}
By taking the limit $q\to 1$, the algebra $\C[N]^{N'(w)}$ is therefore endowed with the structure of a Poisson 
algebra (see \eg \cite[III.5.4]{BG}), and (\ref{eq-q-commut}) yields
\begin{equation}\label{eq_Poisson}
 \{\varphi_{T_i}, \varphi_{T_j}\} = \l_{ij} \varphi_{T_i}\varphi_{T_j}, \qquad (1\le i,j \le r),
\end{equation}
that is, the Poisson structure on $\C[N]^{N'(w)}$ is compatible with its cluster structure in the sense of \cite[\S4.1]{GSV2}.
Moreover, it was shown \cite[Theorem 12.3]{GLS5} that the $\Q(q)$-algebra $\A_{\Q(q)}(\CC_w)$ 
obtained from $\A_{\Q[q^{\pm1}]}(\CC_w)$ by extending coefficients is isomorphic to the 
quantum coordinate ring $A_q(\n(w))$ coming from the Drinfeld-Jimbo theory of quantum enveloping
algebras. Hence this Poisson structure coincides with the standard Poisson
structure on the Schubert cell $C_w$.

\subsection{}\label{subsect-Poisson} Let now $v\le w$. 
The cluster algebra $R(\CC_{v,w})$ is by construction a subalgebra of 
$\C[N]^{N'(w)}$, and its seeds are certain sub-seeds of the seeds of the cluster structure of 
$\C[N]^{N'(w)}$. It follows that $R(\CC_{v,w})$ has a natural quantization as a quantum cluster subalgebra
$\A_{\Q[q^{\pm1}]}(\CC_{v,w})$ of $\A_{\Q[q^{\pm1}]}(\CC_w)$.
Moreover, the functions $\varphi_P$ with $P$ a $\CC_{v,w}$-projective object are quantized into elements $Y_P$
which $q$-commute with every quantum cluster variable of $\A_{\Q[q^{\pm1}]}(\CC_{v,w})$, thus they form an Ore set
$M_{q,v,w}$ in $\A_{\Q[q^{\pm1}]}(\CC_{v,w})$. Hence we can localize $\A_{\Q[q^{\pm1}]}(\CC_{v,w})$ by $M_{q,v,w}$ and we
obtain in this way a quantization of the cluster algebra $R_{v,w}$ of Theorem~\ref{mainThm}~(iii).
Taking again the limit $q\to 1$, we see that $R_{v,w}$ is endowed with the structure of a Poisson algebra
compatible with its cluster structure.

\section{Examples}\label{sect5}

\subsection{}
We first illustrate the results of \S\ref{sect3} with a small example. 
Let $G$ be of type $A_3$, that is, $G=SL(4)$.
We choose for $H$ the subgroup of diagonal matrices in $SL(4)$,
for $B$ the subgroup of upper triangular matrices, so that $N$ is the subgroup of upper unitriangular matrices. 
Take
\[
 v = s_2, \quad w=s_1s_2s_3.
\]
Then the Schubert cell $C_w$ is 3-dimensional, the opposite Schubert cell $C^v$ is 5-dimensional,
and their intersection $\RR_{v,w}$ is 2-dimensional in the 6-dimensional flag variety $X=B^-\backslash SL(4)$. 
In fact, using for example \cite[\S4]{D} or \cite[Proposition 5.2]{MR}, it is easy to see that $\RR_{v,w}\simeq \C^*\times\C^*$.
On the other hand, a straightforward calculation gives 
\[
N_{v,w}=
\left\{y=\left(\begin{matrix}1&0&t&tu\\0&1&0&0\\0&0&1&u\\0&0&0&1\end{matrix}\right) 
 \mid t,u \in \C^*\right\}
\]
and thus $N_{v,w} \simeq \RR_{v,w}$, in agreement with Lemma~\ref{lemNwK}.
The polynomial functions on $N$ whose non-vanishing defines the open subset $O_{v,w}$ are the minors:
\[
\De_{s_2(\varpi_1),\,s_3s_2s_1(\varpi_1)} = \De_{1,4},\quad
\De_{s_2(\varpi_2),\,s_3s_2s_1(\varpi_2)} = \De_{13,14} = \De_{3,4},\quad
\De_{s_2(\varpi_3),\,s_3s_2s_1(\varpi_3)} = \De_{123,124} = \De_{3,4},
\]
where we have written $\De_{I,J}$ for the minor with row set $I$ and column set $J$.
Hence
\[
 O_{v,w}=\left\{x=\left(\begin{matrix}1&a&b&c\\0&1&d&e\\0&0&1&f\\0&0&0&1\end{matrix}\right)\in N 
 \mid c\not = 0,\ f\not = 0\right\}.
\]
An elementary calculation gives, for $x\in O_{v,w}$, 
\[
 \zeta_{v,w}(x) = \left(\begin{matrix}
 1&\ds\frac{-1}{c}&\ds\frac{-a}{c}&\ds\frac{-b}{c}\\[4mm]
 0&1&a&\ds\frac{bf-c}{f}\\[4mm]
 0&0&1&\ds\frac{df-e}{f}\\[4mm]
 0&0&0&1\end{matrix}\right).
\]
One then checks easily that $\zeta_{v,w}(nx)=\zeta_{v,w}(x)$ for $n\in N(v)$, 
hence $\zeta_{v,w}(O_{v,w}) = \zeta_{v,w}(\O_{v,w})$,
where 
\[
\O_{v,w} = O_{v,w}\cap N'(v) =
\left\{x=\left(\begin{matrix}1&a&b&c\\0&1&0&e\\0&0&1&f\\0&0&0&1\end{matrix}\right) 
 \mid c\not = 0, f\not = 0\right\},
\]
and, moreover, the restriction of $\zeta_{v,w}$ to $\Omega_{v,w}$ is injective.
Now, $x\in\O_{v,w}$ satisfies $\zeta_{v,w}(x)\in N(w^{-1})$, where
\[
 N(w^{-1})=\left\{y=\left(\begin{matrix}1&\a&\b&\gamma\\0&1&0&0\\0&0&1&0\\0&0&0&1\end{matrix}\right) 
 \mid \a,\b,\gamma \in \C\right\},
\]
if and only if $x$ is of the form
\[
x =  \left(\begin{matrix}
 1&0&\ds\frac{c}{f}&c\\[4mm]
 0&1&0&0\\[1mm]
 0&0&1&{f}\\[1mm]
 0&0&0&1\end{matrix}\right),
\]
that is, if and only if $x\in N_{v,w}$, in agreement with Lemma~\ref{lemalg} and Lemma~\ref{lempenul}.
Finally, the action of $N(v)\times N'(w)$ on $N$ is given by
\[
nxn'=
\left(\begin{matrix}
 1&0&0&0\\
 0&1&A&0\\
 0&0&1&0\\
 0&0&0&1\end{matrix}\right) 
\left(\begin{matrix}
1&a&b&c\\
0&1&d&e\\
0&0&1&f\\
0&0&0&1\end{matrix}\right)
\left(\begin{matrix}
1&B&C&0\\
0&1&D&0\\
0&0&1&0\\
0&0&0&1\end{matrix}\right) 
=
\left(\begin{matrix}
1&a'&b'&c'\\
0&1&d'&e'\\
0&0&1&f'\\
0&0&0&1\end{matrix}\right)
=x'
 \]
where
\[
a'=a+B,\quad 
b'=b+C+Da,  \quad 
c'=c,\quad 
d'=d+A+D,\quad 
e'=e+Af,\quad 
f'=f.
\]
If $x\in O_{v,w}$ then $c,f\in\C^*$. So if $x'=x$ for $x\in O_{v,w}$
we have 
\[
B=0,\quad
A=0,\quad
D=0,\quad
C=0
\]
which forces $n=n'=e$, in agreement with Theorem~\ref{theo2}.
(Note on the other hand that if $f=0$ the action of $N(v)\times N'(w)$ on $x$ is \emph{not} free.)
Finally, we see that 
\[
{}^{N(v)}\C[O_{v,w}]^{N'(w)} = \C[c^{\pm 1},f^{\pm1}] = \C[N_{v,w}]. 
\]

\subsection{}\label{sect-example}
We now illustrate the results of \S\ref{sect2}, \S\ref{sec-decategorif} and \S\ref{sect-Poisson}.
Let $G$ be of type $A_5$, that is, $G=SL(6)$. We take
\[
 v = s_1s_2s_1s_4s_5s_4, \quad w=s_1s_3s_5s_2s_4s_1s_3s_5s_2s_4s_1s_3s_5s_4 = w_0s_2.
\]
Then the Schubert cell $C_w$ is 14-dimensional, the opposite Schubert cell $C^v$ is 9-dimensional,
and their intersection $\RR_{v,w}$ is 8-dimensional in the 15-dimensional flag variety $X=B^-\backslash SL(6)$.

\subsubsection{}
We consider the Frobenius category
\[
\CC_{v,w} = \CC^v\cap\CC_w=\E^\dag_{v^{-1}}(\md(\L)) \cap \E_u(\md(\L)), 
\]
with $u=w^{-1}w_0 = s_2$. We fix the reduced word 
$\bi = (i_{14},\ldots,i_1) = (1,3,5,2,4,1,3,5,2,4,1,3,5,4)$. 
The corresponding increasing sequence $(v_{(0)},\ldots,v_{(14)})$ defined in \S\ref{defvj} is 
\[
\begin{array}{l}
v_{(0)} = e,\quad
v_{(1)}= s_4,\quad
v_{(2)}= v_{(3)}=s_5s_4,\quad
v_{(4)}= s_1s_5s_4,\quad
v_{(5)}= s_4s_1s_5s_4,\\[2mm]
v_{(6)}= v_{(7)}=v_{(8)}=s_2s_4s_1s_5s_4,\quad
v_{(9)}= v_{(10)}=v_{(11)}=v_{(12)}=v_{(13)}=v_{(14)}=v.
\end{array}
\]
It follows that $J_{v,\bi} = \{3,7,8,10,11,12,13,14\}$.
Then, according to Proposition~\ref{prop-desc-U}, the indecomposable direct summands of 
the $\CC_{v,w}$-cluster-tilting object $U_\bi = V_{\bi}/t_v(V_\bi)$ look as follows:
\begin{align*}
U_3 &= {\bsm&4\\3\esm} &
U_7 &= {\bsm3\esm} &
U_8 &= {\bsm&&&4\\1&&3&&5\\&2&&4\\&&3\esm} &
U_{10} &= {\bsm1&&3\\&2&&4\\&&3\esm}
\\[3mm]
U_{11} &= {\bsm&&4\\&3&&5\\2&&4\\&3\esm} &
U_{12} &= {\bsm1\\&2\\&&&3\esm} &
U_{13} &= {\bsm&&3\\&2&&4\\1&&3&&5\\&2&&4\\&&3\esm} &
U_{14} &= {\bsm&&5\\&4\\3\esm}
\end{align*}
The numbers can be interpreted as basis vectors or as composition factors (see \cite[\S2.4]{GLS4}).
The summands $U_{10}$, $U_{11}$, $U_{12}$, $U_{13}$, $U_{14}$ are the indecomposable 
$\CC_{v,w}$-projective-injective objects.

By Corollary~\ref{cor-number}, the corresponding cluster variables of $R(\CC_{v,w})$ are the 
following minors of an upper unitriangular $6\times 6$ matrix
\begin{align*}
x_3&=\De_{s_4s_5(\varpi_3),\ s_4s_5s_3(\varpi_3)}, &
x_7&=\De_{s_4s_5s_1s_4s_2(\varpi_5),\ s_4s_5s_3s_1s_4s_2s_5(\varpi_5)},
\\
x_8&=\De_{s_4s_5s_1s_4s_2(\varpi_3),\ s_4s_5s_3s_1s_4s_2s_5s_3(\varpi_3)}, &
x_{10}&=\De_{s_4s_5s_1s_4s_2s_1(\varpi_4),\ s_4s_5s_3s_1s_4s_2s_5s_3s_1s_4(\varpi_4)},
\\
x_{11}&=\De_{s_4s_5s_1s_4s_2s_1(\varpi_2),\ s_4s_5s_3s_1s_4s_2s_5s_3s_1s_4s_2(\varpi_2)}, &
x_{12}&=\De_{s_4s_5s_1s_4s_2s_1(\varpi_5),\ s_4s_5s_3s_1s_4s_2s_5s_3s_1s_4s_2s_5(\varpi_5)},
\\
x_{13}&=\De_{s_4s_5s_1s_4s_2s_1(\varpi_3),\ s_4s_5s_3s_1s_4s_2s_5s_3s_1s_4s_2s_5s_3(\varpi_3)}, &
x_{14}&=\De_{s_4s_5s_1s_4s_2s_1(\varpi_1),\ s_4s_5s_3s_1s_4s_2s_5s_3s_1s_4s_2s_5s_3s_1(\varpi_1)}.
\end{align*}
Using more usual notation and writing $\De_{I,J}$ for the minor with row set $I$ and column set $J$,
we have
\begin{align*}
x_3&=\De_{3,\ 5}, &
x_7&=\De_{3,\ 4}, &
x_8&=\De_{123,\ 256}, &
x_{10}&=\De_{123,\ 245},
\\
x_{11}&=\De_{23,\ 56},&
x_{12}&=\De_{123,\ 234},&
x_{13}&=\De_{123,\ 456},&
x_{14}&=\De_{3,\ 6}.
\end{align*}
The cluster variables $x_{10}$, $x_{11}$, $x_{12}$, $x_{13}$, $x_{14}$ are frozen.

To describe the cluster seed containing this cluster, we need to attach a quiver $\G_\bi$ to it.
This quiver is determined from the cluster-tilting module $U_{\bi}$ as follows.
Its vertices are the 8 indecomposable direct summands of $U_\bi$. The number of arrows
from $U_j$ to $U_k$ is the dimension of the space of \emph{irreducible} morphisms from $U_j$
to $U_k$ in the category $\add(U_\bi)$, that is, morphisms that cannot be factored non trivially
in $\add(U_\bi)$.
Here is the quiver $\G_\ii$:
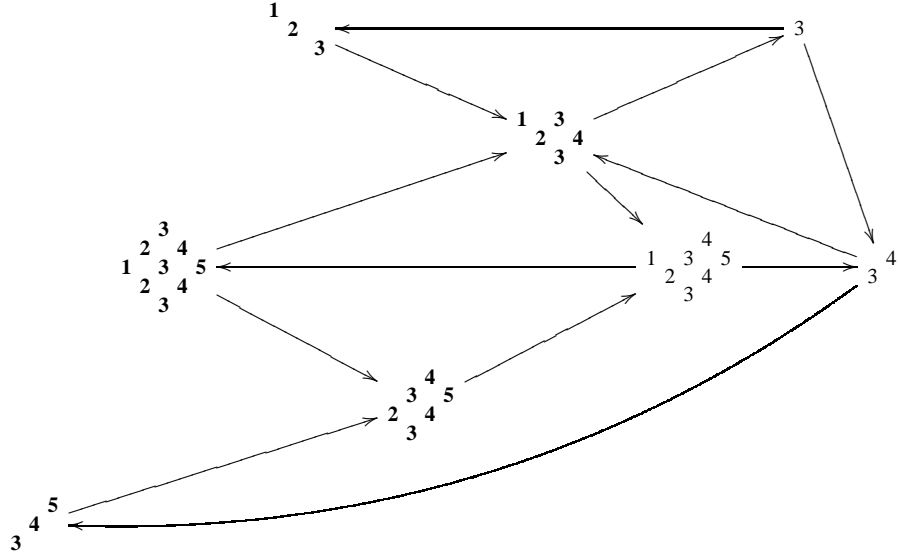
\[
\xymatrix@-0.7pc{
&&{\bsm\bf1\\&\bf2\\&&&\bf3\esm}\ar@{<-}[rrrr]
&&&&{\bsm3\esm}\ar@{<-}[dll]
\\
&&&&{\bsm\bf1&&\bf3\\&\bf2&&\bf4\\&&\bf3\esm}\ar@{<-}[dlll]\ar@{<-}[drrr]\ar@{<-}[ull]
\\
&{\bsm&&\bf3\\&\bf2&&\bf4\\\bf1&&\bf3&&\bf5\\&\bf2&&\bf4\\&&\bf3\esm}\ar@{<-}[rrrr]
&&&&{\bsm&&&4\\1&&3&&5\\&2&&4\\&&3\esm}\ar@{<-}[ul]\ar@{<-}[dll]
&&{\bsm&4\\3\esm}\ar@{<-}[ll]\ar@{<-}[uul]
\\
&&&{\bsm&&\bf4\\&\bf3&&\bf5\\\bf2&&\bf4\\&\bf3\esm}\ar@{<-}[ull]\ar@{<-}[dlll]
\\
{\bsm&&\bf5\\&\bf4\\\bf3\esm}\ar@/_1.1cm/@{<-}[uurrrrrrr]
}
\]
For example $U_7 = {\bsm3\esm}$ maps to the socle of 
$U_{10} = {\bsm1&&3\\&2&&4\\&&3\esm}$, but this map factors through
$U_{12} = {\bsm1\\&2\\&&&3\esm}$, so it is not irreducible.
On the other hand $U_{10}$ maps to $U_7$ with kernel ${\bsm1\\&2&&4\\&&3\esm}$ 
and this map is irreducible. Hence we get one arrow from $U_{10}$ to $U_7$.

The modules printed in bold face are the indecomposable $\CC_{v,w}$-projective-injectives, and therefore
they give frozen vertices of $\G_{\bi}$. If we remove them, as well as all incident arrows,
we are left with the quiver
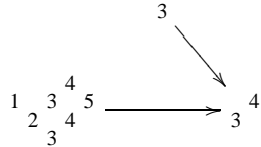
\[
\xymatrix@-0.7pc{
&{\bsm3\esm}
\\
{\bsm&&&4\\1&&3&&5\\&2&&4\\&&3\esm}
&&{\bsm&4\\3\esm}\ar@{<-}[ll]\ar@{<-}[ul]
}
\]
of type $A_3$. Thus, $\CC_{v,w}$ has finitely many indecomposable objects and $R(\CC_{v,w})\subset\C[N]$
is a cluster algebra of finite type $A_3$.
So, by Theorem~\ref{mainThm}~(iv), the coordinate ring $\C[\RR_{v,w}]$ is isomorphic to the 
cluster algebra $R_{v,w}$ obtained by inverting the frozen variables of $R(\CC_{v,w})$.

\subsubsection{}
The cluster variables $\{x_j \mid j\in J_{v,\bi}\}$ form a system of log-canonical coordinates for the Poisson structure
on $\RR_{v,w}$ introduced in \S\ref{subsect-Poisson}. 
Using (\ref{log-can}), (\ref{eq_Poisson}), and the explicit description of the $\L$-modules $U_j\ (j\in J_{v,\bi})$,
we can easily calculate the matrix $L=(\l_{ij})$ such that
\[
 \{x_i,x_j\}= \l_{ij}x_ix_j,\qquad (i,j\in J_{v,\bi}).
\]
We obtain:
\[
L = 
\left(\begin{matrix}
       0&-1&0&0&0&0&0&1\\
       1&0&1&0&1&1&0&1\\
       0&-1&0&-1&-1&-1&0&0\\
       0&0&1&0&0&0&0&1\\
       0&-1&1&0&0&0&0&0\\
       0&-1&1&0&0&0&0&0\\
       0&0&0&0&0&0&0&0\\
       -1&-1&0&-1&0&0&0&0
      \end{matrix}
\right).
\]

\subsubsection{}\label{example-positroid}
We can get a more concrete description of $\C[\RR_{v,w}]$ if we notice that $v=w_K$ is the longest element of the 
maximal parabolic subgroup $W_K$, where $K=\{1,2,4,5\}$. Thus, as noted in \S\ref{sect-partial-flag},
the projection $\pi^K : X\to X_K$ restricts to an isomorphism $\RR_{v,w} \stackrel{\sim}{\to} \RR^K_{v,w}$,
and $R_{v,w}$ is also isomorphic to the coordinate ring of the Lusztig stratum
$\RR^K_{v,w}$. 
Here $X_K = \Gr(3,6)$ is the Grassmannian of $3$-dimensional subspaces of $\C^6$,
and $\RR^K_{v,w}$ is an open positroid variety.
The Grassmannian $X_K$ has a familiar homogeneous coordinate system given by Pl\"ucker coordinates. 
We will denote by $[ijk]$ the Pl\"ucker coordinate labelled by the subset $\{i,j,k\}$
of~$\{1,\ldots,6\}$. Note that $\RR^K_{v,w}\subset \pi^K(C^v)$, and since $v=w_K$ the 
projection $\pi^K(C^v)$ is the open cell of $\Gr(3,6)$ given by the non-vanishing of $[123]$.
We may therefore fix $[123]=1$ and consider the remaining Pl\"ucker coordinates as 
affine coordinates on $\pi^K(C^v)$. It is then easy to check that $\RR^K_{v,w}$ is
the locally closed subset of $\pi^K(C^v)$ defined by 
\begin{equation}\label{eq-plucker}
[345]=0,\quad 
[234]\not = 0,\quad
[245]\not = 0,\quad
[456]\not = 0,\quad
[156]\not = 0,\quad
[126]\not = 0,
\end{equation}
(the non-vanishing of $[123]$ is understood).
Now recall that the coordinate ring of $\pi^K(C^v)$ is isomorphic to ${}^{N(v)}\C[N]$,
by mapping the minor $\De_{123,\ ijk}$ to the Pl\"ucker coordinate $[ijk]$.
We denote by $\tR_{v,w}$ the cluster algebra structure on $\C[\RR^K_{v,w}]$ deduced
from $R_{v,w}$ under this isomorphism.
Notice that the non-vanishing Pl\"ucker coordinates of (\ref{eq-plucker})
exactly correspond to the frozen variables
\[
x_{10}=\De_{123,\ 245},\quad
x_{11}=\De_{23,\ 56} =\De_{123,\ 156},\quad
x_{12}=\De_{123,\ 234},\quad
x_{13}=\De_{123,\ 456},\quad
x_{14}=\De_{3,\ 6} =\De_{123,\ 126}.
\]
The cluster algebra $\tR_{v,w}$ has $9$ (non frozen) cluster variables, 8 of which are the Pl\"ucker
coordinates
\[
[125],\quad
[124],\quad
[256],\quad
[145],\quad
[235],\quad
[146],\quad
[236],\quad
[246].
\]
(The cluster variables $x_3$, $x_7$, $x_8$ of $R_{v,w}$ correspond respectively to $[125]$, $[124]$, $[256]$.)
The 9th cluster variable $x$ of $R_{v,w}$ comes from the unique indecomposable object of $\CC_{v,w}$ with a 2-dimensional
socle, namely
\[
 M = {\bsm1&&3&&5\\&22&&44\\&&33\esm}.
\]
In view of (\ref{eq-cluster-char2}), the two short exact sequences in $\CC_{v,w}$
\[
0\to {\bsm&3&\\2&&4\\&3\esm} \to M \to  {\bsm1&&&&5\\&2&&4\\&&3\esm}\to 0,
\qquad
0\to {\bsm1&&&&5\\&2&&4\\&&3\esm}\to  {\bsm&&3\\&2&&4\\1&&3&&5\\&2&&4\\&&3\esm} \to {\bsm&3&\\2&&4\\&3\esm} \to 0
\]
allow to express $x=\varphi_M$ in terms of minors. 
Using again the isomorphism  
${}^{N(v)}\C[N]\stackrel{\sim}{\to} \C[\pi^K(C^v)]$, we get that $x$ corresponds to
the cluster variable
\[
 \tilde{x} = [145][236]-[456]
\]
of $\tR_{v,w}$.
Finally, note that the restrictions to $\RR^K_{v,w}$ of the remaining 5 Pl\"ucker coordinates:
\[
[134],\qquad
[135],\qquad
[136],\qquad
[346],\qquad
[356],
\]
which are \emph{not} cluster variables in $\tR_{v,w}$, can nevertheless be expressed as Laurent monomials
in the cluster variables, if we use Pl\"ucker relations together with (\ref{eq-plucker}). Indeed
we have in $\C[\RR^K_{v,w}]$:
\begin{align*}
 [134]&=\ds\frac{[145][234]}{[245]},&
 [135]&=\ds\frac{[145][235]}{[245]},&
 [346]&=\ds\frac{[234][456]}{[245]},
 \\[3mm]
 [356]&=\ds\frac{[235][456]}{[245]},&
 [136]&=\ds\frac{\tilde{x}}{[245]}.
\end{align*}
This shows that if the 8 cluster variables of any given cluster evaluate positively at a given point $x\in\RR^K_{v,w}$, 
then all the Pl\"ucker coordinates (except [345]) also have positive evaluation at $x$. 
In other words, the positive part $\RR^K_{v,w}$ in the sense of Postnikov \cite{P} or Lusztig \cite{L} coincides
with the positive part induced by our cluster structure.

\subsection{}\label{subsect-compare}
Inspired by Postnikov \cite{P}, Muller and Speyer have conjectured a cluster algebra structure on the coordinate 
rings of open positroid varieties, described in terms of alternating strand diagrams \cite[Conjecture 3.4]{MS}. 
It is an interesting question to compare
this cluster structure with the cluster structure given by Theorem~\ref{mainThm}.

Let us do it for the particular example developed in \S\ref{example-positroid},
that is, for the open positroid variety $\RR^K_{v,w}$ of the Grassmannian $X_K = \Gr(3,6)$ labelled
by the pair $(v,w)$, where $v=w_K$ and $w=w_0s_2$.
In that case, the cluster algebra $\A$ of \cite[Conjecture 3.4]{MS} is as follows (see the running example of 
\cite[\S3]{MS}).
The frozen cluster variables of $\A$ are: 
\[
f_1 = [234],\quad
f_2 = [346],\quad
f_3 = [456],\quad
f_4 = [156],\quad
f_5 = [126].
\]
(Here we omit, as we may, the frozen variable $f_0 = [123]$, because we prefer to work in the affine chart $[123]=1$,
as in \S\ref{example-positroid}.) The mutable cluster variables are the 8 Pl\"ucker coordinates
\[
[124],\quad
[134],\quad
[136],\quad
[236],\quad
[356],\quad
[246],\quad
[146],\quad
[256],  
\]
together with the degree 2 element
\[
y = [124][356]-[456] = [125][346]. 
\]
An initial seed for $\A$ is given by the cluster $([134],[136],[356])$ and the quiver displayed in
\cite[Figure 4]{MS}.
Thus, we see that the cluster algebra $\A$ is of finite type $A_3$, like the cluster algebra $\tR_{v,w}$ of 
\S\ref{example-positroid}. But we also see that $\tR_{v,w}$ and $\A$ are \emph{different}: they have different sets 
of cluster variables. Note in particular that the frozen variable $f_2=[346]$ of $\A$ is \emph{not} a cluster variable
of $\tR_{v,w}$.
However, using the formulas of \S\ref{example-positroid}, we can express all the cluster variables
of $\A$ as Laurent monomials in the cluster variables of $\tR_{v,w}$, with non-vanishing denominator on $\RR_{v,w}^K$,
and vice-versa. Therefore we get two different cluster algebra structures on the coordinate ring $\C[\RR_{v,w}^K]$,
but they give rise to the \emph{same positive part} of the open positroid variety.

On the other hand, in \cite[Remark 3.5]{MS} Muller and Speyer have proposed a ``mirror image'' of their conjecture, in which
the strands of the Postnikov diagrams are labelled by their \emph{source} instead of their \emph{target}. 
In our running example, this alternative cluster algebra $\A'$ would exactly coincide with our cluster
algebra $\tR_{v,w}$.

\section*{Main notation}

For the convenience of the reader we collect here the list of the main notation of the paper, and we indicate where it is defined.

\begin{itemize}
 \item[\S1.1] $G$, $H$, $B$, $B^-$, $W$, $w_0$, $X$, $\pi$, $C_w$, $C^w$, $\RR_{v,w}$.
 \item[\S1.2] $I$, $x_i(t)$, $y_i(t)$, $B_K^-$, $X_K$, $\pi_K$, $\pi^K$, $W_K$, $W^K$, $w_K$, $\RR_{w_Kv,w}^K$.
 \item[\S1.4] $N$, $\n$.
 \item[\S2.1] $N(w)$, $N'(w)$, $x^\Tr$, $[z]^+_w$, $G_0$. 
 \item[\S2.2] $\bw$, $\bbw$, $\vpi_i$, $\De^{\vpi_i}$, $\De_{u(\vpi_i),v(\vpi_i)}$.
 \item[\S2.3] $N_{v,w}$.
 \item[\S2.4] $O_{v,w}$.
 \item[\S2.5] $U_{v,w}$, $\zeta_{v,w}$.
 \item[\S2.6] $\Omega_{v,w}$.
 \item[\S2.7] $\phi_w$, ${}^{N(v)}\C[O_{v,w}]^{N'(w)}$.
 \item[\S2.8] $D_{v,w}$, $M_{v,w}$, ${}^{N(v)}\C[N]^{N'(w)}$, $S_{v,w}$, 
 \item[\S3.1] $\mod(A)$, $\Sub(M)$, $\Fac(M)$, $\underline{\CC}$.
 \item[\S3.2] $\L$, $S_i$, $Q_i$, $\tau$, $s_i$, $\E_i$, $\E_i^\dag$, $\E_w$, $\E_w^\dag$, $I_w$, $\CC_w$, $J_w$, $\CC^w$,
 $t_w$.
 \item[\S3.3] $\CC_{v,w}$.
 \item[\S3.4] $\add(T)$, $\mu_M(T)$, $\Gamma_T$.
 \item[\S3.5] $I_{i,w}$, $Q_{i,v,w}$.
 \item[\S4.1] $\varphi_M$.
 \item[\S4.4] $\Sigma_T$, $R(\B)$, $S(\B)$.
 \item[\S4.6] $\ZZ$, $\vph_Z$, $\SS^*$.
 \item[\S4.8] $R$, $R^+$, $\gamma_k$, $V_k$, $V_\bi$, $U_\bi$, $R(\CC_{v,w})$, $v_{(k)}$, $J_{v,\bi}$, $U_j$, 
 $\Sigma(M)$, $R_{v,w}$, $\widetilde{R}_{v,w}$.
 \item[\S5.1] $N_v^w$, $\Delta_v^w$.
 \item[\S6.1] $\A_{\Q[q^{\pm1}]}(\CC_w)$, $Y_{T_i}$, $\la_{ij}$.
 \item[\S6.2] $\A_{\Q[q^{\pm1}]}(\CC_{v,w})$, $M_{q,v,w}$.
\end{itemize}

\bigskip
\bigskip
{\bf\large Acknowledgements}

\bigskip
The project of understanding the strata of the nonnegative part of a 
partial flag variety in terms of cluster algebras was started by Nicolas Chevalier in his PhD thesis 
\cite{C2} defended in 2012. He introduced the categories $\CC_{v,w}$ in
the particular case when $v=w_K$ is the maximal element of a parabolic
subgroup of $W$, and he conjectured that $\CC_{w_K,w}$ gives a categorical
model for a cluster algebra structure in the coordinate ring of 
$\RR^K_{w_K,w}$. He checked the conjecture in a number of examples when
$X_K$ is a type $A$ Grassmannian (in particular by comparing it with \cite{BFZ}
when $\RR^K_{w_K,w}$ is isomorphic to a double-Bruhat cell for $GL_n$).
I am very grateful to Nicolas for allowing me
to use his work as a starting point for this paper. I have also benifited
a lot from discussions with Milen Yakimov, who motivated me for pursuing
Nicolas' project. In particular, he suggested that the restriction $v=w_K$ in \cite{C2}
might not be necessary. I want to thank Pierre Baumann, Michel Brion and Milen Yakimov for their
valuable remarks on a preliminary version of this paper.
Finally, I am grateful to Greg Muller and David Speyer for informing me of 
Remark 3.5 in the new version of \cite{MS}.


\bigskip
\small
\noindent
\begin{tabular}{ll}
Bernard {\sc Leclerc} : & Normandie Univ, France;\\ 
&UNICAEN, LMNO F-14032 Caen, France;\\
&CNRS UMR 6139, F-14032 Caen, France;\\
&Institut Universitaire de France.\\
&email : {\tt bernard.leclerc@unicaen.fr}
\end{tabular}
\end{document}